\documentclass[11pt,leqno,english]{article}
 \usepackage{latexsym}
 \usepackage{amssymb}
 \usepackage{amsmath}
 \usepackage{amsthm}
\usepackage{natbib}

 \usepackage[dvips]{graphics}
\usepackage{graphicx}

 \newcommand{ \un }{\mathbb{I} }
 \newcommand{ \p }{\mathbb{P} }
 \newcommand{ \pa }{\mathbb{P}^{\alpha} }
  \newcommand{ \pao }{\mathbb{P}^{\alpha}_{0} }
\newcommand{ \pam }{\mathbb{P}^{\alpha}_{m_n}}
 \newcommand{ \E }{\mathbb{E}}
 \newcommand{ \Ea }{\mathbb{E}^{\alpha}}
 \newcommand{ \Eam }{\mathbb{E}^{\alpha}_{m_n}}

\newcommand{ \F }{ \mathbb{F} }
 
 \newcommand{ \Z }{ \mathbb{Z} }
 \newcommand{\N}{ \mathbb{N} }
 
 \newcommand{ \V }{\textrm{Var} }
 
 \newcommand{ \rw}{\textrm{R.W.R.E. }}

 \newcommand{ \Ct }{ \mathcal{C} }
 
 \newcommand{ \f }{ \mathcal{F} }

 \newcommand{ \tm }{ m }
 \newcommand{ \tM }{ M }
\newcommand{ \tmo }{ m_n }

\newcommand{ \Ie }{ I_{\eta_0} }
\newcommand{ \eto }{ \eta_0 }

\newcommand{ \D}{ \mathcal{D} }
\newcommand{ \A}{ \mathcal{A} }

 \newcommand{ \lo }{ \mathcal{L} }

\newtheorem{The}{{\bf Theorem}}[section]
\theoremstyle{definition}
\newtheorem{Def}[The]{{\bf Definition}}
\theoremstyle{plain}
 \newtheorem{Lem}[The]{Lemma}
 \newtheorem{Cor}[The]{\bf Corollary}
 \newtheorem{Pro}[The]{\bf Proposition}
 \theoremstyle{definition}
\newtheorem{Rem}[The]{{\bf Remark}}

 \newenvironment{Pre}{\noindent \textbf{Proof.} \\ }{$\
 \blacksquare$}

\newenvironment{Preth}[1]{\noindent \textbf{Proof (of Theorem #1).} \\ }{$\
 \blacksquare$}
\newenvironment{Prepr}[1]{\noindent \textbf{Proof (of Proposition #1).} \\ }{$\
 \blacksquare$}

\newenvironment{Prele}[1]{\noindent \textbf{Proof (of Lemma #1).} \\ }{$\
 \blacksquare$}

\makeatletter

\@addtoreset{equation}{section} \makeatother

\title{ On the concentration of Sinai's walk
\\ \vspace{1cm}
 \large{Pierre Andreoletti} $^{\dag}$,\footnote{ Centre de Physique Th\'eorique, C.N.R.S. UMR 6207, Universit\'e Aix-Marseille I, II, Universit\'e du Sud - Toulon - Var, F.R.U.M.A.M. (Marseille, France)
 and Centro de Modelamiento Mathematico C.N.R.S. U.M.R 2071, Universidad de Chile (Santiago,
 Chile). \newline \vspace{0.1cm}  $\quad$  MSC 2000 60G50; 60J55. \newline \vspace{0.5cm} \textit{Key words and phrases :  Random environment, random walk,
Sinai's regime, Markov chain, local time.} \newline \vspace{0.5cm}
CPT-2004/P.068}   }

\begin{document}

\bibliographystyle{unsrtnat}
\maketitle

\noindent  $^{\dag}$ Universit\'e Aix-Marseille II, Facult\'e des
sciences de Luminy, C.P.T. case 907, 13288 Marseille cedex 09
France. e-mail : \texttt{andreole@cpt.univ-mrs.fr}

\noindent \\ \textbf{Abstract:} We consider Sinai's random walk in
random environment. We prove that for an interval of time $[1,n]$
Sinai's walk sojourns in a small neighborhood of the point of
localization for the quasi totality of this amount of time.
Moreover the local time at the point of localization normalized by
$n$ converges in probability to a well defined random variable of
the environment. From these results we get applications to the favorite sites of the walk and to the maximum of the local time.

\section{Introduction}

Random Walks in Random Environment (R.W.R.E.) are basic processes
in random media. The one dimensional case with nearest neighbor
jumps, introduced by \cite{Solomon}, was first studied by
\cite{KesKozSpi}, \cite{Sinai}, \cite{Golosov},  \cite{Golosov0}
and \cite{Kesten2} all these works show the diversity of the
possible behaviors of such walks depending on hypothesis assumed
for the environment. At the end of the eighties \cite{Deh&Revesz}
and \cite{Revesz} give the first almost sure behavior of the \rw
in the recurrent case. Then we have to wait until the middle of
the nineties to see new results. An important part of these new
results concerns the problem of large deviations first studied by
\cite{GreHol2} and then by \cite{ZeiGan}, \cite{PiPO},
\cite{PiPoZe} and \cite{ZeCoGa} (see \cite{Zeitouni} for a
review). In the same period using the stochastic calculus for the
recurrent case \cite{Shi},
 \cite{HuShi2}, \cite{HuShi1}, \cite{Hu1}, \cite{Hu} and \cite{HuShi0} follow the works of
\cite{Schumacher} and \cite{Brox}  to give very precise results on
the random walk and its local time (see \cite{Shi1} for an
introduction). Moreover recent results on the problem of aging are
given in \cite{ZeDEGu}, on the moderate deviations in
\cite{ComPop} for the recurrent case, and on the local time in
\cite{GanShi} for the transient case. In parallel to all these
results a continuous time model has been studied, see for example
\cite{Schumacher} and \cite{Brox}, the works of \cite{Tanaka},
\cite{Mathieu2}, \cite{Tanaka2}, \cite{KawTan}, \cite{Mathieu1}
and \cite{Taleb0}.

Since the beginning of the eighties the delicate case of \rw in
dimension larger than 2 has been studied a lot. For recent reviews
(before 2002) on this topics see the papers of \cite{ Sznitman}
and \cite{Zeitouni}. See also, \cite{Sznitman3}, \cite{Varadhan},
\cite{Rass2} and \cite{ZeiCom2}.

In this paper we are interested in Sinai's walk i.e the one
dimensional random walk in random environment with three
conditions on the random environment: two necessaries hypothesis
to get a recurrent process (see \cite{Solomon}) which is not a
simple random walk and an hypothesis of regularity which allows us
to have a good control on the fluctuations of the random
environment.

The asymptotic behavior of such walk was discovered by
\cite{Sinai}, he showed that this process is sub-diffusive and
that at time $n$ it is localized in the neighborhood of a well
defined point of the lattice. This \textit{point of localization}
 is a random variable depending only on the random environment
and $n$, its explicit limit distribution was given, independently,
by \cite{Kesten2} and A.
 \cite{Golosov0}.

We prove, with a probability very near one that this process is concentrated in a small neighborhood of the point of localization,
this means that for an interval of time $[1,n]$ Sinai's walk
spends the quasi totality of this amount of time in the
neighborhood of the point of localization. The size of this
neighborhood $\thickapprox (\log \log  n)^2$ is negligible
comparing to the typical range $(\log n)^2$ of Sinai's walk. Extending this result to a neighborhood  of arbitrary size we get that, with a strong probability, the size of the interval where the walk spends more than a half of its time is smaller than every positive strictly increasing sequence. We also
 prove that the local time of this random walk at the point of
localization normalized by $n$ converges in probability to a
random variable depending only on $n$ and on the random
environment. This random variable is the inverse of the mean of
the local time at the valley where the walk is trapped within a
return time to the point of localization, we prove that the mean with respect to the environment of this mean is bounded. We generalize this
result for neighboring points of the point of localization. All
our results are "quenched" results, this means that we work with a
fixed environment that belongs to a probability subset of the
random environment that has a probability that goes to one as $n$
 diverges. We give some consequences on the maximum of the local time and
the favorite sites of Sinai's walk.

\noindent \\ This paper is organized as follows. In section 2 we
describe the model and recall Sinai's results. In section
\ref{sec2} we present our main results. In sections \ref{sec4b}
 and \ref{sec4} we give the proof of these results.


\section{Description of the model and Sinai's results}

\subsection{Sinai's random walk definition }

Let $\alpha =(\alpha_i,i\in \Z)$ be a sequence of i.i.d. random
variables taking values in $(0,1)$ defined on the probability
space $(\Omega_1,\f_1,Q)$, this sequence will be called random
environment. A random walk in random environment (denoted
R.W.R.E.) $(X_n,n
\in\N)$ is a sequence of random variable taking value in $\Z$, defined on $( \Omega,\f,\p)$ such that \\
$ \bullet $ for every fixed environment $\alpha$, $(X_n,n\in \N)$
 is a Markov chain with the following transition probabilities, for
 all $n\geq 1$ and $i\in \Z$
 \begin{eqnarray}
& & \p^{\alpha}\left[X_n=i+1|X_{n-1}=i\right]=\alpha_i, \label{mt} \\
& & \p^{\alpha}\left[X_n=i-1|X_{n-1}=i\right]=1-\alpha_i \equiv
\beta_i. \nonumber
\end{eqnarray}
We denote $(\Omega_2,\f_2,\pa)$ the probability space
associated to this Markov chain. \\
 $\bullet$ $\Omega = \Omega_1 \times \Omega_2$, $\forall A_1 \in \f_1$ and $\forall A_2 \in \f_2$,
$\p\left[A_1\times
A_2\right]=\int_{A_1}Q(dw_1)\int_{A_2}\p^{\alpha(w_1)}(dw_2)$.

\noindent \\ The probability measure $\pa\left[\left.
.\right|X_0=a \right]$  will be  denoted $\pa_a\left[.\right]$,
 the expectation associated to $\pa_a$: $\Ea_a$, and the expectation associated to $Q$:
 $\E_Q$.

\noindent \\ Now we introduce the hypothesis we will use in all
this work. The two following hypothesis are the necessaries
hypothesis
\begin{eqnarray}
 \E_Q\left[ \log
\frac{1-\alpha_0}{\alpha_0}\right]=0 , \label{hyp1bb} \label{hyp1}
\end{eqnarray}
\begin{eqnarray}
\V_Q\left[ \log \frac{1-\alpha_0}{\alpha_0}\right]\equiv \sigma^2
>0 . \label{hyp0}
\end{eqnarray}
 \cite{Solomon} shows that under \ref{hyp1} the
process $(X_n,n\in \N)$ is $\p$ almost surely recurrent and
\ref{hyp0} implies that the model is not reduced to the simple
random walk. In addition to \ref{hyp1} and \ref{hyp0} we will
consider the following hypothesis of regularity, there exists $0<
\eta_0 < 1/2$ such that
\begin{eqnarray}
& & \sup \left\{x,\ Q\left[\alpha_0 \geq x \right]=1\right\}= \sup
\left\{x,\ Q\left[\alpha_0 \leq 1-x \right]=1\right\} \geq \eta_0.
\label{hyp4}
\end{eqnarray}

\noindent We call \textit{Sinai's random walk} the random walk in
random environment previously defined with the three hypothesis
\ref{hyp1}, \ref{hyp0} and \ref{hyp4}.

\subsection{ The random potential and the valleys }

Let
\begin{eqnarray}
\epsilon_i \equiv \log \frac{1-\alpha_i}{\alpha_i},\ i\in \Z,
\end{eqnarray}
define :
\begin{Def} \label{defpot2} The random potential $(S_m,\  m \in
\Z)$ associated to the random environment $\alpha$ is defined in the following way: for all $k$ and $j$, if $k>j$
\begin{eqnarray} 
 &&S_k-S_{j}=\left\{ \begin{array}{ll} \sum_{j+1\leq i \leq k} \epsilon_i, &  k\neq0, \\
  - \sum_{j \leq i \leq -1} \epsilon_i , &  k=0  , \end{array} \right. \nonumber \\
&& S_{0}=0, \nonumber
 \end{eqnarray}
and symmetrically if $k<j$.
\end{Def}

\begin{figure}[h]
\begin{center}
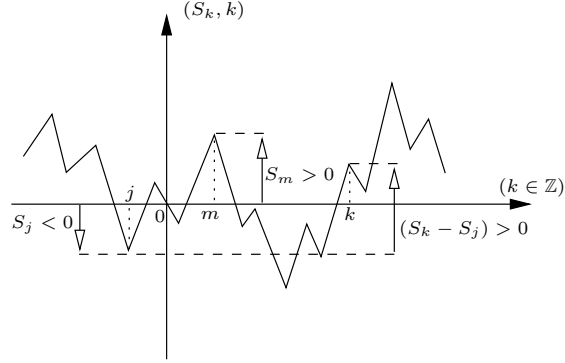 \caption{Trajectory of the random potential} \label{fig5}
\end{center}
\end{figure}

\noindent
\begin{Rem}  
using Definition \ref{defpot2} we have :
\begin{eqnarray} 
 S_k=\left\{ \begin{array}{ll} \sum_{1\leq i \leq k} \epsilon_i, &  k=1,2,\cdots , \\
   \sum_{k \leq i \leq -1} \epsilon_i , &  k=-1,-2,\cdots  , \end{array}  \right. \label{defmerde}
 \end{eqnarray}
however, if we use \ref{defmerde} for the definition of $(S_{k},k)$, $\epsilon_{0}$ does not appear in this definition and moreover it is not clear, when $j<0<k$, what the difference $S_k-S_{j}$ means (see figure \ref{fig5}).   
\end{Rem}

\begin{Def} \label{c2s2d1}
 We will say that the triplet $\{M',m,M''\}$ is a \textit{valley} if
 \begin{eqnarray}
& &  S_{M'}=\max_{M' \leq t \leq m} S_t ,  \\
& &  S_{M''}=\max_{m \leq t \leq
 \tilde{M''}}S_t ,\\
& & S_{m}=\min_{M' \leq t \leq M''}S_t \ \label{2eq58}.
 \end{eqnarray}
If $m$ is not unique we choose the one with the smallest absolute
value.
 \end{Def}

\begin{Def} \label{deprofvalb}
 We will call \textit{depth of the valley} $\{\tM',\tm,\tM''\}$ and we
 will denote it $d([M',M''])$ the quantity
\begin{eqnarray}
 \min(S_{M'}-S_{m},S_{M''}-S_{m})
 .
 \end{eqnarray}
 \end{Def}

\noindent  Now we define the operation of \textit{refinement}
 \begin{Def}
Let  $\{\tM',\tm,\tM''\}$ be a valley and let
  $\tM_1$ and $\tm_1$ be such that $\tm \leq \tM_1< \tm_1 \leq \tM''$
  and
 \begin{eqnarray}
 S_{\tM_1}-S_{\tm_1}=\max_{\tm \leq t' \leq t'' \leq
 \tM''}(S_{t'}-S_{t''}) .
 \end{eqnarray}
 We say that the couple $(\tm_1,\tM_1)$ is obtained by a \textit{right refinement} of $\{\tM',\tm,\tM''\}$. If the couple $(\tm_1,\tM_1)$ is not
 unique, we will take  the one such that $\tm_1$ and $\tM_1$ have the smallest  absolute value. In a similar way we
  define the \textit{left refinement} operation.
 \end{Def}

\begin{figure}[h]
\begin{center}
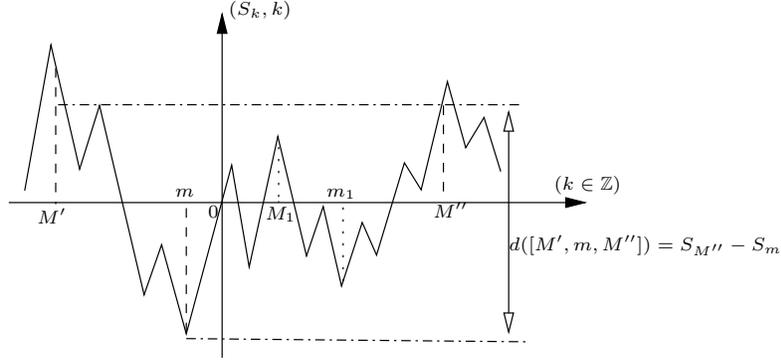 \caption{Depth of a valley and refinement operation} \label{fig4}
\end{center}
\end{figure}

\noindent \\  We denote $\log_p $ with  p $\geq 2$ the $p$
iterated logarithm. In all this work we will suppose that $n$ is
large enough such that $\log_p n$ is positive.

\begin{Def} \label{thdefval1b} For $\gamma>0$, $n>3$ and $\Gamma_n \equiv \log n+ \gamma \log_2 n $, we say that a valley $\{\tM',\tm,\tM''\}$ contains $0$ and
is of depth larger than $\Gamma_n$ if and only if
\begin{enumerate}
\item  $ 0 \in [\tM',\tM'']$, \item $d\left([\tM',\tM'']\right)
\geq \Gamma_n $ , \item if $\tm<0,\ S_{\tM''}-\max_{ \tm \leq t
\leq 0
}\left(S_t\right) \geq \gamma \log_2 n $ , \\
 if $\tm>0,\
S_{\tM'}-\max_{0 \leq t \leq
 \tm}\left(S_t\right) \geq \gamma \log_2 n $ .
\end{enumerate} $\gamma$ is a free parameter.
\end{Def}

\subsection{The basic valley $\{{M_n}',\tmo,{M_n}\}$}

 We recall the notion of \textit{ basic valley } introduced
by Sinai and denoted here $\{{M_n}',\tmo,{M_n}\}$. The definition
we give is inspired by the work of \cite{Kesten2}. First let
$\{\tM',\tmo,\tM''\}$ be the smallest \textit{valley that contains
$0$ and of depth larger than} $\Gamma_n$. Here smallest means that
if we construct, with the operation of refinement, other valleys
in $\{\tM',\tmo,\tM''\}$ such valleys will not satisfy one of  the
properties of Definition \ref{thdefval1b}. ${M_n}'$ and ${M_n}$
are defined from $\tmo$ in the following way: if $\tmo>0$
\begin{eqnarray}
& & {M_n}'=\sup \left\{l\in \Z_-,\ l<\tmo,\ S_l-S_{\tmo}\geq
\Gamma_n,\
S_{l}-\max_{0 \leq k \leq \tmo}S_k \geq \gamma \log_2 n \right\} ,\\
& & {M_n}=\inf \left\{l\in \Z_+,\ l>\tmo,\ S_l-S_{\tmo}\geq
\Gamma_n\right\} . \label{4.8}
\end{eqnarray}
if $\tmo<0$
\begin{eqnarray}
& & {M_n}'=\sup \left\{l\in \Z_-,\ l<\tmo,\ S_l-S_{\tmo}\geq
\Gamma_n\right\} , \\
& & {M_n}=\inf \left\{l\in \Z_+,\ l>\tmo,\ S_l-S_{\tmo}\geq
\Gamma_n,\ S_{l}-\max_{ \tmo \leq k \leq 0}S_k \geq \gamma \log_2
n \right\} . \label{4.10}
\end{eqnarray}
if $\tmo=0$
\begin{eqnarray}
& & {M_n}'=\sup \left\{l\in \Z_-,\ l<0,\ S_l-S_{\tmo}\geq
\Gamma_n \right\} , \\
& &  {M_n}=\inf \left\{l\in \Z_+,\ l>0,\ S_l-S_{\tmo}\geq \Gamma_n
\right\} . \label{4.12}
\end{eqnarray}
\noindent  $\{{M_n}',\tmo,{M_n}\}$ exists
 with a $Q$ probability as close to one as we need. In fact it is not
 difficult to prove the following lemma (see Section \ref{2par22} for the ideas of the proof)

\begin{figure}[h]
\begin{center}
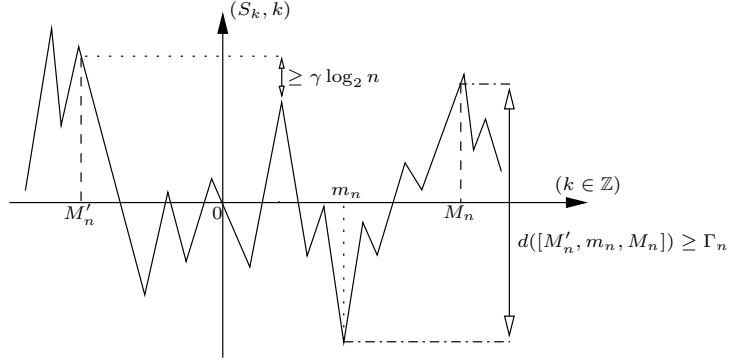 \caption{Basic valley, case $m_{n}>0$ } \label{thfig4}
\end{center}
\end{figure}

\begin{Lem} \label{moexiste} There exists $c>0$ such that if  \ref{hyp1}, \ref{hyp0} and \ref{hyp4} hold, for all $\gamma>0$ there exists $n_0\equiv n_0(\gamma,Q)$
such that for all $n>n_0$
\begin{eqnarray}
Q\left[\{{M_n}',\tmo,{M_n}\} \neq \varnothing \right] \geq
1-\frac{c \gamma \log_2 n}{\log n}.
\end{eqnarray}
\end{Lem}

\noindent In all this paper we use the same notation $n_0$ for an
integer that could change from line to line. Moreover in the rest
of the paper we do not always make explicit the dependence on the
free parameter $\gamma$ and on the distribution $Q$ of all those $n_0$ even if Lemma
\ref{moexiste} is constantly used.

\subsection{ Localization phenomena (\cite{Sinai})}
\begin{The} \label{thSinai1}
Assume \ref{hyp1}, \ref{hyp0} and \ref{hyp4} hold. For all
$\gamma>6$, $\epsilon>0$ and $\delta>0$ there exists $n_0$ such
that for all $n>n_0$, there exists $G_n \subset \Omega_1$ with
$Q\left[G_n\right] \geq 1- \epsilon $ such that
\begin{eqnarray}
\lim_{n\rightarrow +\infty } \sup_{\alpha \in G_n}
\pa_0\left[\left|X_n -m_n\right|
> \delta (\log n)^2 \right]=0
\end{eqnarray}
\end{The}

\noindent This result shows that with a $Q$ and $\pa$ probability
 as close to 1 as we want, at time $n$ the \rw is localized in a
small neighborhood of $\tmo $. The parameter $\gamma$ comes from
the definition of $\tmo \equiv \tmo(\gamma)$.

\noindent  Y. Sinai shows also, with a probability close to one,
that for a given time interval $[0,n]$, the \rw is trapped in the
basic valley and therefore is sub-diffusive. In fact if we define
\begin{eqnarray}
W_n=\left\{{M_n}',\ {M_n}'+1,\cdots,\ \tmo,\cdots,\ {M_n}-1,\
{M_n}\right\},
\end{eqnarray}
we have
\begin{Pro}\label{lem2} Assume \ref{hyp1}, \ref{hyp0} and \ref{hyp4} hold. For all
$\gamma>6$, $\epsilon>0$ there exists $E\equiv E(\epsilon)$ and
$n_0$ such that for all $n>n_0$, there exists $G_n \subset
\Omega_1$ with $Q\left[G_n\right] \geq 1- \epsilon $ such that
 \begin{eqnarray}
 \lim_{n\rightarrow +\infty
} \inf_{\alpha \in G_n}
\p^{\alpha}_0\left[\bigcap_{m=0}^{n}\left\{X_m\in W_n
\right\}\right] =1 ,
 \end{eqnarray}
moreover
\begin{eqnarray}
 \lim_{n\rightarrow +\infty
} \inf_{\alpha \in G_n}
\p^{\alpha}_0\left[\bigcap_{m=0}^{n}\left\{X_m\in\left[-E(
\sigma^{-1} \log n)^2 , E( \sigma^{-1} \log n)^2
\right]\right\}\right] =1 .
 \end{eqnarray}
\end{Pro}

\noindent See for example \cite{Zeitouni} or \cite{Pierre1} for
alternative proofs of these results.

\section{Main results: concentration phenomena \label{sec2}}
Let us define the local time at $k$ $(k\in \Z)$ within the
interval of time $[1,T]$ ($T \in \N^*$) of $(X_n,n\in \N)$
\begin{eqnarray}
\lo\left(k,T\right) \equiv \sum_{i=1}^T \un_{\{X_i=k\}} .
\end{eqnarray}
$\un$ is the indicator function ($k$ and $T$ can be deterministic or random variables). Let $V\subset \Z$, we denote
\begin{eqnarray}
\lo\left(V,T\right) \equiv \sum_{j \in V} \lo\left(j,T\right)
=\sum_{i=1}^T\sum_{j \in V} \un_{\{X_i=j\}} .
\end{eqnarray}

\subsection{Local time in a neighborhood of $\tmo$}

Let us define the following sequences, let $p\geq 2$ and $n$ large enough, define : :
\begin{eqnarray}
&& f_{p}(n)=[(\log_2 n \log_p n)^2]\ ([-] \textrm{ is the integer
part of } -) \label{F}, \\
& & R_{p}(n)= (\log_{p+1} n )^{ 1/2} (\log_{p} n )^{-1/2}.
\end{eqnarray} 
For all $n \geq 1$ define the set :
\begin{eqnarray}
 & & \F_{p}(n)=\{\tmo-f_{p}(n),\tmo-f_{p}(n)+1,\
\cdots,\tmo,\tmo+1,\ \cdots,\
 \tmo+f_{p}(n)\}. \label{V}
\end{eqnarray}

\begin{The} \label{tpslocmo} There exists $c>0$ such that if
  \ref{hyp1bb}, \ref{hyp0} and \ref{hyp4} hold,
 for all $p\geq 2$, $0< \rho < 2$ and $\gamma \geq 11$, there exists $c_1 \equiv c_1(Q,\gamma)>0$  and $n_0$ such that
for all $n>n_0$ there exists $G_n' \subset \Omega_1$ with
$Q\left[G_n'\right] \geq 1-c  R_{p}(n)
 $ and
\begin{eqnarray}
\inf_{\alpha \in
G_n'}\left\{\pao\left[\frac{\lo\left(\F_{p}(n),n\right)}{n} \geq
\left(1-\frac{1}{(f_{p}(n))^{\rho/2}}\right)\right]
\right\}\geq 1-\frac{c_1}{ (f_{p}(n) )^{1-\rho/2}}.
\end{eqnarray}
\end{The}

\begin{Rem}  \label{remrembb}
The set $G_n'$ called the \textit{set of good environments} will
be defined in section \ref{sec4.5}. \\
\noindent In fact $n_0$ and $c_{1}$ depends only on $Q$ through
 $\Ie \equiv \log [(1-\eta_0)/(\eta_0)],
Q\left[\epsilon_0>\Ie/2\right]$,
 $Q\left[\epsilon_0<-\Ie/2\right]$, $\sigma$ and $\E_Q[|\epsilon_0|^3]$, however to simplify the writing we do not make it explicit. In all this work we denote
$c$ a strictly positive numerical constant that can grow from line
to line
 if needed.\\
If we choose $p=2$ in Theorem \ref{tpslocmo} we get better rates
for the convergence of the probabilities, however using \ref{V} we
get that $|\F_{p}(n)| = (\log_2 n)^4$ ($|A|$ denotes the cardinal of
 $A \subset \Z$), whereas for $p>2$, $|\F_{p}(n)|= (\log_2 n)^2
(\log_{p} n)^2$. Recall that $\gamma$ comes from the definition of
the point $\tmo \equiv \tmo(\gamma)$, the condition $\gamma \geq
11$ will become clear in Section \ref{condgamma} when we will
prove this theorem.
\end{Rem}


\noindent In addition to Y. Sinai's results on localization, with
a $\pa$ and $Q$ probability very near one, in a time interval
$[0,n]$ the \rw does not spend a finite proportion of this time
interval outside $\F_{p}(n)$. In fact the time spend outside $\F_{p}(n)$ does
not exceed $n/(\log_2 n \log_p n)^{\rho}$. Moreover, notice that $
|\F_{p}(n)| /|W_n| \approx (\log_2 n \log_p n)^2 (\log n)^{-2} \searrow
0$ that is the subset $\F_{p}(n)$ of $W_n$ where the \rw stays a time
greater than $n(1-1/(\log_2 n \log_p n)^{\rho})$ has no density
inside $W_n$ in the limit when $n$ goes to
infinity. \\

\begin{Rem} Notice that if we look for an annealed result, that means a result in $\p$ probability, with the same condition of Theorem \ref{tpslocmo} we get
\begin{eqnarray}
\p\left[\frac{\lo\left(\F_{p}(n),n\right)}{n} \geq
\left(1-\frac{1}{(f_{p}(n))^{\rho/2}}\right)\right]
\geq 1-\frac{c_1}{ (f_{p}(n) )^{1-\rho/2}}-cR_{p}(n).
\end{eqnarray}
By definition $\lo\left(\F_{p}(n),n\right) \leq 1$, therefore we also have
\begin{eqnarray}
\p\left[\left| \frac{\lo\left(\F_{p}(n),n\right)}{n}-1\right| \leq
\frac{1}{(f_{p}(n))^{\rho/2}}\right]
\geq 1-\frac{c_1}{ (f_{p}(n) )^{1-\rho/2}}-cR_{p}(n).
\end{eqnarray}
\end{Rem}

\noindent \\ To prove this theorem we will prove the following
 key result on the environment, define
\begin{eqnarray}
& &  T_{x}=\left\{\begin{array}{l} \inf\{k\in\N^*,\ X_k=x \}
\\ + \infty \textrm{, if such }  k \textrm{ does not exist.}
 \end{array} \right. \label{3.7sec}
\end{eqnarray}
\begin{Pro} \label{propoquidonneb} \label{propoquidonne} There exists $c>0$ such that if \ref{hyp1bb}, \ref{hyp0} and
\ref{hyp4} hold, for all $p\geq 2$ there exists $n_0$ such that
for all $n>n_0$
\begin{eqnarray}
Q\left[\Ea_{\tmo}\left[ \lo(\bar{\F}_{p}(n),T_{\tmo})\right]>
\frac{2}{\eto}\frac{1}{f_{p}(n)+1}\right] \leq c R_{p}(n) \label{2.7}
\end{eqnarray}
where $\bar{\F}_{p}(n)$ is the complementary of $\F_{p}(n)$ in $W_{n}$.
\end{Pro}

\noindent This implies that with a $Q$ probability going to one
when $n$ goes to infinity, in $\Ea$ mean, the \rw will never reach
a point $k \in \bar{\F}_{p}(n)$ in a return time to $\tmo$ in spite of
the that $ |\bar{\F}_{p}(n)| /|W_n| \approx 1$. The proof of this
proposition is given Section \ref{5.2}.

\noindent \\
If we want to replace in Theorem  \ref{tpslocmo} the neighborhood $\F_{p}(n)$ by an \textbf{arbitrary small}  neighborhood $(\Theta_{n},n)$  but such that $\lo(\Theta_{n},n)/n$ converge in probability to
one, we get with a similar method the following result:
\begin{The} \label{tpslocmos} There exists $c>0$ such that if
  \ref{hyp1bb}, \ref{hyp0} and \ref{hyp4} hold,
 for all $0<\rho<1/4$, all strictly positive increasing sequences $(\theta(n),n\geq 1)$ and $\gamma \geq 11$, there exists $c_{1}\equiv c_{1}(Q,\gamma)$ and $ n_0 $ such that
for all $n>n_0$ there exists, $G_n' \subset \Omega_1$ with
$Q\left[G_n'\right] \geq 1-c(\theta(n))^{-\rho}
 $ and
\begin{eqnarray}
\inf_{\alpha \in
G_n'}\left\{\pao\left[ \lo(\Theta(n),n) \geq n\left(1-\frac{1}{(\theta(n))^{\rho}}\right)  \right]
\right\}\geq 1- \frac{c_{1}}{(\theta(n))^{1/2-2\rho}}, \label{3.9}
\end{eqnarray}
where $\Theta(n)=\{m_{n}-\theta(n),m_{n}-\theta(n)+1, \  \cdots, m_{n},m_{n+1},\cdots,m_{n}+\theta(n)\}$.
\end{The}

\noindent The key result on the environment for the proof of this result is the following, assume that $(\theta(n),n)$ is a strictly positive increasing sequence, 
\begin{Pro} \label{propoquidonnebbb}  Assume \ref{hyp1bb}, \ref{hyp0} and
\ref{hyp4} hold, for all $0<\rho<1/2$ there exists $c_{1}\equiv c_{1}(Q)$ and $n_0$ such that
for all $n>n_0$
\begin{eqnarray}
Q\left[\Ea_{\tmo}\left[ \lo(\bar{\Theta}(n),T_{\tmo})\right]>
\frac{1}{(\theta(n))^{1/2-\rho}}\right] \leq \frac{
c_{1}}{ (\theta(n))^{\rho} } \label{2.7b}
\end{eqnarray}
where $\bar{\Theta}(n)$ is the complementary of $\Theta(n)$ in $W_{n}$.
\end{Pro}

\noindent We call this facts (Theorem \ref{tpslocmo} and
\ref{tpslocmos}) \textit{concentration} even if it could be more
appropriate to define the concentration in term of the random
variable
\begin{eqnarray}
Y_n=\inf_{x \in \Z}\min\left\{k>0\ :\ \lo([x-k,x+k],n)>n/2
\right\}.
\end{eqnarray}
In fact Theorem \ref{tpslocmos} implies the following result
\begin{The} \label{tpslocmoss} There exists $c>0$  such that if
  \ref{hyp1bb}, \ref{hyp0} and \ref{hyp4} hold,
 for all strictly positive increasing sequences $(\theta(n),n)$, there exists $c_{1}\equiv c_{1}(Q)>0$ and $ n_0 $ such that
for all $n>n_0$ there exists $G_n' \subset \Omega_1$ with
$Q\left[G_n'\right] \geq 1-c(\theta(n))^{-1/4}
 $ and
\begin{eqnarray}
\inf_{\alpha \in G_n'}\left\{\pao\left[Y_{n} \leq   \theta(n)
\right] \right\}\geq 1-
\frac{c_{1}}{(\theta(n))^{1/4}}.\end{eqnarray}
\end{The}

\noindent 
Notice that the annealed result is given by
\begin{eqnarray}
\p\left[Y_{n} \leq   \theta(n)
\right]  \geq 1-
\frac{c_{2}}{(\theta(n))^{1/4}},
\end{eqnarray}
with $c_{2}\equiv c_{2}(Q)=c+c_{1}$.

\noindent  It would be now interesting to study the $\p$
almost sure behavior of $Y_n$, these questions are beyond the scope of
this paper, see \cite{Pierre3} for a first approach.


\subsection{Local time on $\tmo$ and on its neighboring points}
\noindent \\ The following theorem is a result of weak law of
large number type for the local time of the \rw at $\tmo$. We
obtain the convergence of the local time at $\tmo$ normalized by
$n$ to a
 $Q$ random variable as Y. Sinai obtain the convergence of $X_n$
 to $m_n$ which is also a $Q$ random variable.

\begin{The} \label{tpslocmo1} There exists $c>0$ such that if \ref{hyp1bb}, \ref{hyp0} and \ref{hyp4} hold, for all $\gamma \geq 11$,
there exists $c_{1}\equiv c_{1}(Q,\gamma)$ and $n_0$ such that for all $n>n_0$ there
exists $G_n' \subset \Omega_1$ with $Q\left[G_n'\right] \geq 1-cR_{2}(n)$ and
\begin{eqnarray}
 & \quad &  \sup_{\alpha \in G_n'}
\left\{\pao\left[\left|\frac{\lo(\tmo,n)}{n}-
\frac{1}{\Eam\left[\lo(W_n,T_{\tmo})\right]} \right|> \frac{1}{
(\log_2 n)^{2}} \right] \right\} \leq \frac{c_1}{ (\log_2
n)^{2}}. \label{10eq43}
\end{eqnarray}
\end{The}



 \noindent \\ We can easily give an intuitive idea of this
result. From Proposition \ref{lem2} we know that the \rw, within
an  interval of time $[1,n]$, spends all its time in the valley
$\{{M_n}',\tmo,{M_n}\}$ with a probability very near one. So, in
spite of the fact that $\Eam\left[T_{\tmo}\right]=+\infty \
Q.a.s.$ as it is easy to check,
 until the instant $n$, the mean of the return time to $\tmo$ is, heuristically, of order $\Eam\left[\lo(W_n,T_{\tmo})\right]$. So we can chop the interval $[1,n]$ in
$\lo(\tmo,n)$ pieces of length
$\Eam\left[\lo(W_n,T_{\tmo})\right]$, therefore
\begin{eqnarray*}
& & n \thickapprox \lo(\tmo,n) \Eam\left[\lo(W_n,T_{\tmo})\right]
\Leftrightarrow \frac{\lo(\tmo,n)}{n} \thickapprox
\left(\Eam\left[\lo(W_n,T_{\tmo})\right]\right)^{-1} .
\end{eqnarray*}

\noindent Now let us give some precisions on the random variable
$\Eam\left[\lo(W_n,T_{\tmo})\right]$, first for all fixed
environment we have the following explicit formula :
\begin{Pro} \label{Pro} For all fixed $\alpha$ and $n$, we have
\begin{eqnarray*}
\Eam\left[\lo(W_n,T_{\tmo})\right]=1&+& \sum_{k=\tmo+1}^{M_n}
\frac{\alpha_{\tmo}}{\beta_k}
\frac{\sum_{j=\tmo+1}^{k-1}\exp(S_{j}-S_{k})}{\sum_{j=\tmo+1}^{k-1}\exp(S_{j}-S_{\tmo})}
\\ &+ & \sum^{\tmo-1}_{k=M_n'} \frac{\beta_{\tmo}}{\alpha_k}
\frac{\sum_{j=k-1}^{\tmo+1}\exp(S_{j}-S_{k})}{\sum_{j=k-1}^{\tmo+1}\exp(S_{j}-S_{\tmo})}.
\end{eqnarray*}
\end{Pro}

\begin{Rem}
Trivially we have $\Eam\left[\lo(W_n,T_{\tmo})\right] \geq 1$ because $\tmo \in W_n$ and
$\lo(\tmo,T_{\tmo})=1$.
\end{Rem}

\noindent We see that  $\Eam\left[\lo(W_n,T_{\tmo})\right] $ depends on the random
environment
 in a complicated way, however using the hypothesis \ref{hyp4} we easily 
 prove the following result :
\begin{Pro} \label{pro3.10} With a $Q$ probability equal to one, for all $n$ we have :
\begin{eqnarray*}
 \frac{\eta_0}{1-\eta_0}\sum_{k \in W_n,k \neq\tmo  }
\frac{1}{\exp(S_{k}-S_{\tmo})} \leq
\Eam\left[\lo(W_n,T_{\tmo})\right]-1 \leq \frac{1}{\eta_0}\sum_{k
\in W_n,k \neq\tmo } \frac{1}{ \exp(S_{k}-S_{\tmo})}.
\end{eqnarray*}
\end{Pro}
\noindent  So we can have a good idea of the fluctuations of
$\Eam\left[\lo(W_n,T_{\tmo})\right]$ only by studying the random variable
\begin{eqnarray}
\sum_{k \in W_n,k \neq\tmo } \frac{1}{ \exp(S_{k}-S_{\tmo})}.
\label{rev1}
\end{eqnarray}

\noindent  The following result is a key result for the random
environment :
\begin{Pro} \label{3pro61b} There exists $c_1 \equiv c_1(Q) >0$ such that
for all $n \geq 1$
\begin{eqnarray}
1 \leq \E_Q\left[ \Eam\left[\lo(W_n,T_{\tmo})\right]\right] \leq
c_1.
\end{eqnarray}
In particular, using Markov inequality, we get that for all
 positive increasing sequences $(\theta(n),n)$  :
\begin{eqnarray}
Q\left[\Ea_{\tmo}\left[ \lo(W_n,T_{\tmo})\right]> \theta(n) \right]
\leq \frac{c_1}{\theta(n)}. \label{eqq}
\end{eqnarray}
\end{Pro}

\begin{Rem}
 $c_1$ depends only on $Q$ through $\eto$ and $Q\left[\epsilon_0<-\Ie/2\right]$ but for simplicity we do not make explicit this dependence, see Section \ref{5.3} for the proof of Proposition
 \ref{3pro61b}.\\
 \noindent \\
One can ask ourself about the convergence in law of
 $\Eam\left[\lo(W_n,T_{\tmo})\right]$ this problem is not an easy consequence of our computations
  and can be by itself an
 independent work.

\end{Rem}

\noindent  One can notice the specificity of the result
\ref{10eq43}, however easy modifications of our method give the
following generalization.

\begin{The} \label{tpslocmo1nn} There exists $c>0$ such that if \ref{hyp1bb},
\ref{hyp0} and \ref{hyp4} hold, for all $\gamma
\geq 11$, there exists $c_{1}\equiv c_{1}(Q,\gamma)$ and $n_0$ such that for all $n>n_0$ there
exists $G_n' \subset \Omega_1$ with $Q\left[G_n'\right] \geq 1-cR_{2}(n) $ and
\begin{eqnarray}
 & & \quad  \sup_{\alpha \in G_n'}
\left\{\pao\left[\bigcup_{k\in
\mathbb{L}(n)}\left\{\left|\frac{\lo(k,n)}{n}-
\frac{1}{\Ea_{k}\left[\lo(W_n,T_{k})\right]} \right|> \frac{ 1 }{
(\log_2 n)^{2}} \right\} \right] \right\} \leq \frac{ c_{1}
\log_{3}n}{ \log_2 n } \label{extension}
\end{eqnarray}
where $\mathbb{L}(n)=\{m_{n}-l(n),m_{n}-l(n)+1, \  \cdots, m_{n},m_{n+1},\cdots,m_{n}+l(n)\}$, $l(n)=\log_{3}(n)/\Ie$ and $\Ie$ is given Remark
\ref{remrembb}.
\end{The}

\begin{Rem}
We deduce this results from the computations we made to prove
Theorem \ref{tpslocmo1}, the point is that, under the hypothesis
of Propositions \ref{propoquidonne} and \ref{3pro61b}, we have the
following results, for all $k\in \mathbb{L}(n)$ :
\begin{eqnarray}
& & Q\left[\Ea_{k}\left[ \lo(\bar{\F}_{p}(n),T_k)\right]>
\frac{2}{\eto}\frac{\log_2 n}{f_{p}(n)+1}\right] \leq c R_{p}(n), \label{lok}
\\
& & Q\left[\Ea_{k}\left[ \lo(\F_{p}(n),T_k)\right]> c_1 (\log_2 n)
\log_{p+1} n  \right] \leq cR_{p}(n). \label{lokk}
\end{eqnarray}
However we think that Theorem \ref{tpslocmo1nn}  is certainly true
for a neighborhood $\mathcal{V}_{n}$ of $\tmo$ $(\F_{p}(n) \subset
\mathcal{V}_{n})$ of size larger than $(\log \log n)^{a}$ for some
$a>0$, but it does not seem to be a simple extension of our
computations.
\end{Rem}

 \noindent The following corollary
 is a simple consequence of \ref{eqq} and Theorem \ref{tpslocmo1}

\begin{Cor}\label{tpslocmo2} There exists
$c>0$ such if \ref{hyp1bb}, \ref{hyp0} and \ref{hyp4} hold,
 for all $p\geq 2$, $\gamma \geq 11$, there exists $c_{1}\equiv c_{1}(Q)$ and $n_0$ such that
for all $n>n_0$ there exists $c_2 \equiv c_{2}(Q) >0$ and  $G_n' \subset \Omega_1$
with $Q\left[G_n'\right] \geq 1-cR_{p}(n) $ and
\begin{eqnarray}
& & \inf_{\alpha \in G_n'} \left\{\pao\left[\lo(\tmo,n) \geq
\frac{n}{2 c_2 \log_{p+1} n} \right] \right\} \geq
1-\frac{c_{1}}{ \log_2 n \log_p
n}, \\
& & \inf_{\alpha \in G_n'} \left\{\pao\left[ \bigcap_{k \in \F_{p}(n)}
\left\{\lo(k,n) \right\}\geq \frac{n}{2 c_2 (\log_2 n) \log_{p+1}
n} \right] \right\} \geq 1-\frac{c_{1} }{ \log_p
n} .
\end{eqnarray}
\end{Cor}

\noindent \\ All these results show that in addition to be
\textit{localized} the \rw is \textit{concentrated} and the region
of concentration and localization are extremely linked together.

\noindent In the following subsection we give some simple
consequences on the maximum of the local time and on the favorite
sites of Sinai's walk.


\subsection{Simple consequences on the maximum of the local time and on the favorite  site of the \rw}

\noindent Let us introduce the following random variables
\begin{eqnarray}
& & \lo^*(n)=\max_{k\in \Z}\left(\lo(k,n)\right) ,\
\tilde{\F}(n)=\left\{k \in \Z,\ \lo(k,n)=\lo^*(n) \right\} .
\end{eqnarray}
$\tilde{\F}(n)$ is the set of all the \textit{favorite sites} and $
\lo^*(n)$ is the maximum of the local times (for a given instant
$n$). The following Corollary is a simple consequence (just inspection) of
\ref{eqq},  Theorem \ref{tpslocmo} and \ref{tpslocmo1}.

\begin{Cor} \label{corder}
There exists $c>0$ such if \ref{hyp1bb}, \ref{hyp0} and \ref{hyp4}
hold,
 for all $p\geq 2$, and $\gamma \geq 11$, there exists $c_{1}\equiv c_{1}(Q,\gamma)$ and $n_0$ such that
for all $n>n_0$ there exists $c_2\equiv c_{2}(Q)>0$ and  $G_n' \subset \Omega_1$
with $Q\left[G_n'\right] \geq 1-cR_{p}(n) $ and
\begin{eqnarray}
& & \inf_{\alpha \in G_n'} \left\{\pao\left[\lo^*(n) \geq
\frac{n}{2c_2 \log_{p+1} n} \right] \right\} \geq
1-\frac{c_{1} }{ \log_2 n \log_p n } , \label{corder1} \\
& & \inf_{\alpha \in G_n'} \left\{ \pao\left[ \tilde{\F}(n) \subset
\F_{p}(n) \right] \right\} \geq 1-\frac{c_{1} }{ \log_2 n
\log_p n}.
\end{eqnarray}
\end{Cor}

\noindent 
To finish we give an interesting application (just inspection) of \ref{corder1}, Theorem \ref{tpslocmos} and \ref{tpslocmo1nn} :

\begin{Cor} \label{Corendend} There exists $c>0$ such that if \ref{hyp1bb},
\ref{hyp0} and \ref{hyp4} hold, for all $\gamma \geq
11$, and all strictly positive increasing sequences $(\theta(n),n)$ there
exists $n_0$ and $c_1 \equiv c_1(Q)$ such that for all $n>n_0$
there exists $G_n' \subset \Omega_1$ with $Q\left[G_n'\right] \geq
1-c (\theta(n))^{-1/4} $ and
\begin{eqnarray}
 & & \quad  \sup_{\alpha \in G_n'}
\left\{\pao\left[\left|\frac{\lo^*(n)}{n}- \max_{k \in
\Theta(n)} \left\{ \frac{1}{\Ea_{k}\left[\lo(W_n,T_{k})\right]}
\right\} \right|> \frac{ 1 }{ (\log_2 n)^{2}} \right]
\right\} \leq \frac{c_1 \theta(n)}{ \log_2 n },
\label{extension}
\end{eqnarray}
recall that $\Theta(n)=\{m_{n}-\theta(n),m_{n}-\theta(n)+1, \  \cdots, m_{n},m_{n+1},\cdots,m_{n}+\theta(n)\}$.
\end{Cor}

\noindent Notice that if we are interested in an annealed result Corollary \ref{Corendend} implies :
\begin{eqnarray}
 & & \p\left[\left|\frac{\lo^*(n)}{n}- \max_{k \in
\Theta(n)} \left\{ \frac{1}{\Ea_{k}\left[\lo(W_n,T_{k})\right]}
\right\} \right|> \frac{ 1 }{ (\log_2 n)^{2}} \right]
 \leq \frac{c_1 \theta(n)}{ \log_2 n }+\frac{c}{(\theta(n))^{1/4}} .
\label{extension2}
\end{eqnarray}

\noindent
The delicate problem of the limit distribution of $\lo^{*}(n)/n$ considered by \cite{Revesz} can not be directly deduced from these results, however they are a good starting point for further investigations on this topic. We think that this limit distribution is very linked to the $Q$ limit distribution of the random variable $ \Eam\left[\lo(W_n,T_{\tmo})\right]$, recall moreover that we have good knowledges on $W_{n}$ thanks to the works of \cite{Kesten2} and \cite{Golosov0}.






\section{Proof of the main results \label{sec4b}}

First we define what we call a \textit{good  environment} and the
\textit{set of good environments.}

\subsection{ Good properties and set of good environments
\label{sec4.5} }

\begin{Def} \label{superb}  Let
$p\geq 2$, $\gamma>0$, $c_1>0$
 and $\omega \in \Omega_1$, we will say that $\alpha \equiv
\alpha(\omega)$ is a \textit{good environment} if there exists
$n_0$ such that for all $n\geq n_0$ the sequence $(\alpha_i,\ i
\in \Z)=(\alpha_i(\omega),\ i \in \Z)$ satisfies
 the properties \ref{3eq325} to \ref{8eq36}
\begin{eqnarray}
& \bullet & \textrm{The valley } \{{M_n}',\tmo, {M_n}\} \textrm{ exists, in particular: } \label{3eq325}\\
& & 0 \in [{M_n}',{M_n}], \\
& & \textrm{ If } \tmo>0,\ S_{{M_n}'}-\max_{0 \leq m \leq \tmo }
\left(S_m\right)\geq \gamma \log_2 n , \\
& & \textrm{ if } \tmo<0,\ S_{{M_n}}-\max_{\tmo \leq m \leq 0 }
\left(S_m\right)\geq \gamma \log_2 n , \\
& &  S_{\tM_n'}-S_{\tmo} \geq \log n+\gamma \log_2 n  ,\label{3eq319}\\
& &  S_{\tM_n}-S_{\tmo} \geq \log n+ \gamma \log_2 n .
\label{3eq320} \\
 & \bullet  & {M_n}'\geq ( \sigma^{-1} \log n)^2
\log_p n ,\ {M_n}\leq ( \sigma^{-1} \log n)^2 \log_p n .
\label{interminibb}
\end{eqnarray}
Define  $\tM_1'$ and $\tm_1'$, respectively the maximizer and
minimizer obtained by the first \textit{left refinement} of the
valley $\{{M_n}',\tmo, {M_n}\}$ and in the same way $\tM_1$ and
$\tm_1$, respectively, the maximizer and minimizer obtained by the
first \textit{right refinement} of the valley $\{{M_n}',\tmo,
{M_n}\}$.
\begin{eqnarray}
& \bullet &  S_{\tM_{1}'}-S_{\tm_1'} \leq \log n-\gamma \log_2 n, \label{3eq326} \\
& &  S_{\tM_{1}}-S_{\tm_1} \leq \log n-\gamma \log_2 n  .
\label{3eq327}
\end{eqnarray}
Define $\F_{p}^{+}(n) = \{\tmo+1,\cdots ,\ \tmo+f_{p}(n)\} ,\ \F_{p}^{-}(n) =
\{\tmo-f_{p}(n),\cdots,\tmo-1\}$ where $f_{p}(n)$ is given by \ref{F}
\begin{eqnarray}
 & \bullet & \min_{k\in
\F_{p}(n)^+}\left(\beta_k
\pa_{k-1}\left[T_k>T_{\tmo}\right]\right) \geq
(g_1(n))^{-1}, \label{8eq38} \\
& & \min_{k\in \F_{p}(n)^-}\left(\alpha_k
\pa_{k+1}\left[T_k>T_{\tmo}\right]\right) \geq
(g_1(n))^{-1}, \label{8eq39}
\end{eqnarray}
where $g_1(n)= \exp\left( \left(
 (4 \sqrt{3} \sigma f_{p}(n) )^2 \log_3 (n) \right )^{1/2} \right) $.
\begin{eqnarray}
 & \bullet & \Ea_{\tmo}\left[ \lo(\F_{p}(n),T_{\tmo}
)\right] \leq c_1 \log_{p+1} n . \label{8eq37} \\
& \bullet & \Ea_{\tmo}\left[ \lo(\bar{\F}_{p}(n),T_{\tmo})\right] \leq
2(\eto(f_{p}(n)+1))^{-1}  . \label{8eq36}
\end{eqnarray}
See  \ref{V} for the definition of $\F_{p}(n)$, and we recall that $\bar{\F}_{p}(n)$ is the complementary of $\F_{p}(n)$ in $W_{n}$.
\end{Def}


\noindent  Define the \textit{set of good environments}
\begin{eqnarray}
G_n'=\left\{\omega \in \Omega_1,\ \alpha(\omega) \textrm{ is a }
\textit{ good environment} \right\}.
\end{eqnarray}
$G_n'$ depends on $p$, $\gamma$, $c_1$ and $n$, however  we do not
make explicit its $p$, $\gamma$ and $c_1$ dependence.

\begin{Pro} \label{profondab} There exists $c>0$ such that if \ref{hyp1bb}, \ref{hyp0} and \ref{hyp4} hold, there exists $c_1>0$ such that for all $p\geq 2$ and $\gamma>0$, there exists $ n_0 \equiv n_0\left(\gamma \right)$
such that for all $n>n_0$
\begin{eqnarray}
Q\left[ G_n'\right] \geq 1-c R_{p}(n)
 .
\end{eqnarray}
\end{Pro}

\begin{Pre}
\noindent  The main ideas of the proof is the subject of section
\ref{sec4}.
\end{Pre}

\noindent \\ For completeness, we recall some results of
\cite{Chung} and \cite{Revesz} on inhomogeneous discrete time
birth and death processes.

\subsection{Basic results on the random walk in a fixed environment}

We will always assume that $\alpha$ is fixed (denoted $\alpha \in
\Omega_1$ in this work). \noindent Let $x \in \Z$, assume $a<x<b$,
the two following lemmata can be found in \cite{Chung} (pages
73-76), the proof follows from the method of difference equations.

\begin{Lem} \label{3.7bb} Recalling \ref{3.7sec}, we have
 \begin{eqnarray}
 & &
 \p^{\alpha}_x\left[T_a>T_b\right]=\frac{\sum_{i=a+1}^{x-1}\exp
\big(S_{i}-  S_{a}\big)  +1}{\sum_{i=a+1}^{b-1}\exp\big( S_{i}-
S_{a}  \big)+1 } \label{k1} , \\  & & \p^{\alpha}_x \left[T_a<T_b
\right]=\frac{\sum_{i=x+1}^{b-1} \exp \big( S_{i}- S_{b} \big)
+1}{\sum_{i=a+1}^{b-1}\exp \big( S_{i}- S_{b} \big) +1 }
\label{k2} .
\end{eqnarray}
\end{Lem}
\noindent Let $T_a \wedge T_b $ be the minimum between $T_a$
and $T_b$.
\begin{Lem} \label{lA0} We have
\begin{eqnarray}
&& \Ea_{a+1}\left[T_{a} \wedge T_{b} \right]=\frac{\sum_{l=a+1}^{b-1}\sum_{j=l}^{b-1}\frac{1}{\alpha_l}F(j,l)}{\sum_{j=a+1}^{b-1}F(j,a)+1}\label{c4s1l1e1} , \\
&& \Ea_x \left[ T_{a} \wedge
T_{b}\right]=\E^{\alpha}_{a+1}\left[T_{a} \wedge T_{b}
\right]\left(1+\sum_{j=a+1}^{x-1}F(j,a)\right)-\sum_{l=a+1}^{x-1}\sum_{j=l}^{x-1}\frac{1}{\alpha_l}F(j,l),
\label{c4s1l1}
\end{eqnarray}
where  $F(j,l)=\exp \big( S_{j}- S_{l} \big) $.
\end{Lem}

\noindent Now we give some explicit expressions for the local
times that can be found in \cite{Revesz} (page 279)

\begin{Lem} \label{4.9}
Under $\pa_x$, $\lo(x,T_{b}\wedge T_{a})$ is a geometric random
variable with parameter
 \begin{eqnarray}
 p=\alpha_{x}\p_{x+1}^{\alpha}\left[T_{x}<T_{b}\right]+\beta_{x}\p^{\alpha}_{x-1}\left[T_{x}<T_{a}\right]
 ,
\end{eqnarray}
that is for all $l \geq 0$, $\pa_{x}\left[\lo(x,T_{b}\wedge
T_{a})=l\right]=p^l(1-p)$.
\end{Lem}

\begin{Lem} \label{Lou} For all $i \in \Z$, we have, if $x>i$
\begin{eqnarray}
\Ea_i\left[
\lo(x,T_i)\right]=\frac{\alpha_i\p_{i+1}^{\alpha}\left[T_{x}<T_{i}\right]}{
\beta_{x}\p^{\alpha}_{x-1}\left[T_{x}>T_{i}\right]} ,
\label{2eq65}
\end{eqnarray}
if $x<i$
\begin{eqnarray}
\Ea_i\left[ \lo(x,T_i)\right]= \frac
{\beta_i\p_{i-1}^{\alpha}\left[T_{x}<T_{i}\right]}{
\alpha_{x}\p^{\alpha}_{x+1}\left[T_{x}>T_{i}\right]} .
\label{2eq68}
\end{eqnarray}
\end{Lem}


\subsection{Proof of Theorem \ref{tpslocmo}, \ref{tpslocmos} and \ref{tpslocmoss} \label{condgamma}}

First we recall the two following elementary results

\begin{Pro}\label{lem2b} There exists $c>0$ such that if \ref{hyp1bb}, \ref{hyp0} and \ref{hyp4}
hold, for all $p\geq 2$ and $\gamma>2$ there exists $n_0 \equiv
n_0\left(\gamma\right)$ such that for all $n>n_0$,
$Q\left[G_n'\right] \geq 1-c R_{p}(n)
 $ and for all $\alpha \in G_n'$
\begin{eqnarray}
\p^{\alpha}_0\left[\bigcup_{m=0}^{n}\left\{X_m\notin W_n
\right\}\right] \leq
 \frac{2 \log_p n}{\sigma^2(\log n)^{\gamma-2}} .  \label{eded}
 \end{eqnarray}
 \end{Pro}

\begin{Pre}
 We use the properties \ref{interminibb}-\ref{3eq327} and  make
 computations
 similar to the ones done by \cite{Sinai}.  The constraint
$\gamma>2$ is here to get a useful result.
\end{Pre}

\begin{Lem}  \label{lemtpsatmo} There exists $c>0$ such that if \ref{hyp1bb}, \ref{hyp0} and \ref{hyp4}
hold, for all $p\geq 2$ and $\gamma \geq 11 $, there exists $n_0
\equiv n_0\left(\gamma\right)$ such that for all $n>n_0$,
$Q\left[G_n'\right] \geq 1-c R_{p}(n)$ and for all $\alpha \in G_n'$
\begin{eqnarray}
 \p^{\alpha}_0\left[ T_{\tmo}> \frac{n}{(\log n)^4} \right]\leq
 \frac{2(\log_p n)^3}{\eto \sigma^6 (\log n)^{\gamma-10}} \label{4.27}.
\end{eqnarray}
\end{Lem}

\begin{Rem}
The constraint $\gamma \geq 11$ that appears in the Theorems
\ref{tpslocmo}, \ref{tpslocmos}, \ref{tpslocmo1}, \ref{tpslocmo1nn} and in their Corollaries comes from
\ref{4.27}.
\end{Rem}

\begin{Pre}
Assume $\tmo>0$, first we remark that
\begin{eqnarray}
\p^{\alpha}_0\left[ T_{\tmo}> \frac{n}{(\log n)^4} \right] \leq
\p^{\alpha}_0\left[ T_{\tmo}> T_{M_n'-1} \right]+
\p^{\alpha}_0\left[ T_{\tmo} \wedge T_{M_n'-1}
> \frac{n}{(\log n)^4} \right] \label{6.3}
\end{eqnarray}
Using \ref{interminibb} and \ref{3eq326} and making computations
similar to the ones done by \cite{Sinai} we get that, for all
$\gamma
> 0$, there exists $n_0 \equiv n_0\left(\gamma\right)$ such that
for all $n>n_0$ and $\alpha \in G_n'$
\begin{eqnarray}
 \p^{\alpha}_0\left[ T_{\tmo}>
T_{M_n'-1} \right] \leq \frac{\log_p n}{ \sigma^2 (\log
n)^{\gamma}} .  \label{6.4}
\end{eqnarray}
Using Markov inequality, Lemma \ref{lA0}, the properties
\ref{3eq326} and \ref{interminibb}
 we obtain for all $\gamma \geq 11$, $n>n_0$ and $\alpha \in G_n'$
\begin{eqnarray}
\p^{\alpha}_0\left[ T_{\tmo} \wedge T_{M_n'-1}
> \frac{n}{(\log n)^4} \right] \leq \frac{(\log_p n)^3}{\eto
\sigma^{6}(\log n)^{\gamma-10}} . \label{6.5}
\end{eqnarray}
Collecting \ref{6.5}, \ref{6.4}, and \ref{6.3} ends the proof of
the  lemma.
\end{Pre}

\begin{Preth}{\ref{tpslocmo}}
 Let $\gamma \geq 11$ and $p \geq 2$, using
Proposition \ref{profondab} we
 take $c>0$ and $n_1$ such that for all  $n\geq n_1$, $Q\left[G_n'\right]
\geq 1-c R_{p}(n)  $. Let
$0<\rho<2$, denote $\delta_n=(f_{p}(n))^{-\rho/2}$, to
prove Theorem \ref{tpslocmo}   we need to give an upper bound of
the
 probability $\pao\left[\lo\left(\F_{p}^c(n),n\right) \geq \delta_n n \right]
$ where $\F_{p}^c(n)$ is the complementary of $\F_{p}(n)$ in $\Z$.
\\ First we use Proposition \ref{lem2b} to
reduce the set $\F_{p}^c(n)$ : there exists $n_2 $ such that for all
$n>n_2$ and $\alpha \in G_n'$
\begin{eqnarray}
& & \pao\left[ \lo\left(\F_{p}^c(n),n\right) > \delta_n n
 \right]  \leq \pao\left[ \lo\left(\bar{\F}_{p}(n),n\right) \geq
\delta_n n \right]+\frac{2 \log_p n}{\sigma^2(\log n)^{\gamma-2}}
\label{4eq8}.
\end{eqnarray}
recall that $\bar{\F}_{p}(n)$ is the complementary of $\F_{p}(n)$ in $W_{n}$. By Lemma
\ref{lemtpsatmo}, there exists $n_3 $ such that for all $n>n_3$
and $\alpha \in G_n'$
\begin{eqnarray}
\quad \pao\left[\lo\left(\bar{\F}_{p}(n),n\right) \geq \delta_n n
\right]  & \leq & \pao\left[\lo\left(\bar{\F}_{p}(n),n\right) \geq
\delta_n n,\ T_{\tmo} \leq \frac{n}{(\log n)^4}  \right]+ \frac{2
(\log_p n)^6}{\eto \sigma^6 (\log n)^{\gamma-10}} \label{4eq11}.
\end{eqnarray}
Let us denote $N_0=\left[n (\log n)^{-4}\right] +1$ and
$\delta_n'=\delta_n-N_0/n$. By the Markov property and the
homogeneity of the Markov chain we get
\begin{eqnarray}
\qquad \pao\left[\lo\left(\bar{\F}_{p}(n),n\right) \geq \delta_n n,\
T_{\tmo} \leq \frac{n}{(\log n)^4}  \right] & \leq &
 \pam\left[ \sum_{k=1}^{n}\un_{\left\{X_k \in
\bar{\F}_{p}(n)\right\} } \geq \delta_n' n
 \right] \label{4eq21} .
\end{eqnarray}
\noindent Let $j\geq 2 $, define the following return times
\begin{eqnarray*}
&& T_{\tmo,j } \equiv \left\{\begin{array}{l}  \inf\{k>T_{\tmo,j-1},\ X_k=\tmo \}, \\
 + \infty \textrm{, if such }  k \textrm{ does not exist.}
 \end{array} \right. \\
& & T_{\tmo,1}  \equiv T_{\tmo} \ (\textrm{see \ref{3.7sec}}).
\end{eqnarray*}
Since by definition $T_{ \tmo,n} > n $, we have
\begin{eqnarray}
 \pam\left[ \sum_{k=1}^{n}\un_{\left\{X_k \in
\bar{\F}_{p}(n)\right\} } \geq   \delta_n' n
 \right]
 & \leq & \pam\left[ \sum_{k=1}^{T_{
\tmo,n}}\un_{\left\{X_k \in \bar{\F}_{p}(n)\right\} } \geq \delta_n' n
\right] . \label{4eq24}
\end{eqnarray}
By definition of the local time and the Markov inequality we get
\begin{eqnarray}
& & \pam\left[ \sum_{k=1}^{T_{\tmo},n}\un_{\left\{X_k \in
\bar{\F}_{p}(n)\right\} } \geq
\delta_n' n \right] \nonumber  \\
& \leq & \left( \sum_{s_1=\tmo+f_{p}(n)+1}^{{M_n}}\Eam\left[
\lo\left(s_1,T_{\tmo,n}\right)\right]+\sum_{s_2={M_n}'}^{\tmo-f_{p}(n)-1}\Eam\left[\lo\left(s_2,T_{\tmo,n}\right)\right]
\right)  (\delta_n' n )^{-1} . \label{5.28}
\end{eqnarray}
Now we use the fact that, by the strong Markov property, the
random variables $\lo\left(s,T_{ \tmo,i+1}-T_{\tmo,i}\right)$ $(0
\leq i \leq n-1)$ are $i.i.d.$, therefore the right hand side of
\ref{5.28} is equal to
\begin{eqnarray}
 \Eam\left[\lo\left(\bar{\F}_{p}(n),T_{\tmo}\right)\right] (\delta_n')^{-1} .\label{4eq35}
\end{eqnarray}
Using  the property \ref{8eq36}, for all  $n>n_1$ and all $\alpha
\in G_n'$
\begin{eqnarray}
\Eam\left[\lo\left(\bar{\F}_{p}(n),T_{\tmo}\right)\right] \leq
\frac{2}{\eto}\frac{1}{f_{p}(n)+1} ,\label{4eq36}
\end{eqnarray}
recall $f_{p}(n)=[(\log_2 n \log_p n)^2]$. Collecting what we did
above, we obtain for all $n \geq n_1 \vee n_2 \vee n_3$ ($a \vee b
=\max(a,b)$)
\begin{eqnarray}
 \pao\left[ \lo\left(\F_{p}^c(n),n\right) > \delta_n n
 \right]
&\leq & \frac{2}{\eto}\frac{1}{(\log_2 n \log_p
n)^2\delta_n'}+\frac{2 \log_p n}{\sigma^2(\log
n)^{\gamma-2}}+\frac{2 (\log_p n)^3}{\eto\sigma^6(\log
n)^{\gamma-10}} . \label{4eq37}
\end{eqnarray}
Taking $\gamma \geq 11$ and choosing $n_0 \geq n_1 \vee n_2 \vee
n_3$ ends the proof of the theorem.\end{Preth}

\noindent \\ In a similar way, we can prove the following
proposition that will be used in the next section

\begin{Pro} \label{thtpsloc1b} There exists $c>0$ such that if \ref{hyp1bb}, \ref{hyp0} and \ref{hyp4}
hold, for all $p\geq 2$, $0<\rho<1$ and $0<\xi'\leq 1$, there
exists $ n_0 $ such that for all $n>n_0$, $Q\left[G_n'\right] \geq
1-c R_{p}(n)
 $ and for all $\alpha \in G_n'$
\begin{eqnarray}
\pam\left[\lo\left(\F_{p}(n), [\xi'n]\right) \geq [n
\xi']\left(1-(\log_2 n (\log_{p}
n)^{2-\rho})^{-1}\right)\right]\geq 1-\frac{2}{\eto}\frac{ 1}{
 (\log_2 n) (\log_p n)^{\rho}} \label{4.40}.
\end{eqnarray}
\end{Pro}

\begin{Preth}{\ref{tpslocmos} and \ref{tpslocmoss}}
First, to prove Theorem \ref{tpslocmos} we use exactly the same method of the proof of Theorem \ref{tpslocmo} but instead of using  Proposition \ref{propoquidonneb} we use Proposition \ref{propoquidonnebbb}.
Now, to get Theorem \ref{tpslocmoss} it is enough to find a point $x$ such that $\lo([x-\theta(n),x+\theta(n)],n) \geq n/2$ in probability, so choosing $x=m_{n}$ and using Theorem \ref{tpslocmos} we get the  result. 
\end{Preth}

\subsection{Proof of Theorem \ref{tpslocmo1} \label{sec5.2}}

\noindent Let $p\geq2$, $\gamma \geq 11$  in all this proof we
take $c>0$, $c_1$ and $n_1\equiv n_1\left(\gamma\right)$ such that
for all $n>n_1$, $Q\left[G_n'\right] \geq 1-cR_{p}(n) $. Let $0<\rho<1$, denote $\nu_n=(\log_{2}n
)^{-1}(\log_p n)^{-\rho} $, $x_1=
\Eam\left[\lo(\F_{p}(n),T_{\tmo})\right]$ and $x_2=
\Eam\left[\lo(\bar{\F}_{p}(n),T_{\tmo})\right]$, by definition we have
\begin{eqnarray}
\Eam\left[\lo(W_n,T_{\tmo})\right]=x_1+x_2. \label{eqeq12b}
\end{eqnarray}
Notice that $\tmo \in \F_{p}(n)$ so $x_1 \geq 1$. Moreover the property
\ref{8eq36} implies that for all $\alpha \in G_n'$ $ x_2 \leq
2(\eta_{0}(f_{p}(n)+1))^{-1}$. We have $
(x_1+x_2)^{-1}=(x_1)^{-1}-x_2(x_1(x_1+x_2))^{-1}$ and one can
 choose $n_2$ such that for all $n>n_2$, $x_2(x_1(x_1+x_2))^{-1}
\leq \nu_n/2 $.  Therefore, for all  $n> n_1 \vee n_2$ and all
$\alpha \in G_n'$, we get
\begin{eqnarray}
\pao \left[\left|\frac{\lo(\tmo,n)}{n}- \frac{1}{x_1+x_2} \right|
> \nu_n \right] & \leq &
 \pao\left[\left|\frac{\lo(\tmo,n)}{n}- \frac{1}{x_1} \right|>
\frac{\nu_n}{2}  \right]. \label{eqeq1}
\end{eqnarray}
 We are left to estimate the right hand side of \ref{eqeq1}, we
 have
\begin{eqnarray}
\qquad \pao\left[\left|\frac{\lo(\tmo,n)}{n}- \frac{1}{x_1}
\right|>\nu_n/2 \right] & \leq & \pao\left[ \lo(\tmo,n)<n \eta_1
\right] + \pao\left[\lo(\tmo,n)>n \eta_2 \right].
\label{3Th1.3.10eq2} \label{3Th1.3.10eq1}
\end{eqnarray}
where $\eta_1 \equiv \eta_1(\rho)=(x_1)^{-1}- \nu_n/2, \
\eta_2\equiv \eta_2(\rho)=(x_1)^{-1}+ \nu_n/2 $.

\noindent \\
 We give an estimate for each terms in the right hand
side of \ref{3Th1.3.10eq1} in the following proposition.

\begin{Pro} \label{8eq33}  There exists $ n_1' $ such
that for all $n>n_1'$ and $\alpha \in G_n'$
\begin{eqnarray} & & \pao\left[\lo(\tmo,n)<\eta_1 n \right]\leq
\frac{4}{\eto}
\frac{ 1}{ \log_2 n (\log_p n)^{\rho}}, \label{5.31} \\
& & \pao\left[\lo(\tmo,n)>\eta_2 n \right]\leq \frac{32(\log_p
n)^3 }{ \eto^2 \sigma^6 (\log n)^{\frac{1}{2}} }. \label{5.32}
\end{eqnarray}
\end{Pro}
\begin{Pre}
We will only prove \ref{5.31}, the proof of \ref{5.32} is easier,
one can check it with a similar method. By Lemma \ref{lemtpsatmo},
there exists
 $n_2
$ such that for all $n>n_2$ and $\alpha \in G_n'$
\begin{eqnarray}
& & \pao\left[\lo\left(\tmo,n\right) < \eta_1n \right]  \leq
\pao\left[\lo\left(\tmo,n\right) < \eta_1 n,\ T_{ \tmo }\leq
\frac{n}{(\log n)^4} \right] +\frac{2(\log_p(n))^3}{\eto \sigma^6
(\log n)^{\gamma-10}} .\label{10eq104}
\end{eqnarray}
Using the strong Markov property and the fact that
$\sum_{j=1}^{T_{\tmo}}\un_{X_j=\tmo}=1$ we get that
\begin{eqnarray}
 \pao\left[\lo\left(\tmo,n\right) < \eta_1 n,\ T_{ \tmo
}\leq \frac{n}{(\log n)^4} \right]  & \leq &
\pam\left[\lo\left(\tmo,(1-\zeta_n)n\right) < \eta_1 n \right]  ,
\label{10eq108}
\end{eqnarray}
where $\zeta_n=N_0/n$, with $N_0=\left[n(\log n)^{-4}\right]+1$.
 Using Proposition \ref{thtpsloc1b} with $\xi'=1-\xi_n$, there exists $n_3 $ such that for all $n>n_3$ and all $\alpha \in G_n'$
\begin{eqnarray}
\qquad &  & \pam\left[\lo\left(\tmo,(1-\zeta_n)n\right) < \eta_1 n
\right] \nonumber
\\
\qquad &   \leq &  \pam\left[\lo\left(\tmo,(1-\zeta_n)n\right) <
\eta_1 n,\ \lo(\F_{p}(n),(1-\zeta_{n})n)\geq (1-\delta_n'')(1-\zeta_n)n
\right] + \frac{ 2(\eto)^{-1}  }{
 (\log_2 n) (\log_p n)^{\rho}} \label{8eq109} .
\end{eqnarray}
where $\delta_n''=(\log_2 n \log_{p+1})^{-1}$. Let us denote
$\eta_1'\equiv \eta_1 ((1-\delta_n'')(1-\zeta_n))^{-1}$, we have
\begin{eqnarray}
\qquad \left\{
\begin{array}{l}
\lo\left(\tmo,n(1-\zeta_n)\right) < \eta_1 n, \textrm{ and} \\
\lo(\F_{p}(n),n(1-\zeta_{n}))\geq (1-\delta_n'')(1-\zeta_n)n.
\end{array} \right.
 & \Rightarrow &
  \lo\left(\tmo,n(1-\zeta_n)\right) <
\eta_1' \lo(\F_{p}(n),n(1-\zeta_{n})) , \label{eq450}
\end{eqnarray}
therefore
\begin{eqnarray}
& & \pam\left[\lo\left(\tmo,n(1-\zeta_n)\right) < \eta_1 n,\
\lo(\F_{p}(n),(1-\zeta_{n})n)\geq (1-\delta_n'')(1-\zeta_n)n
\right] \nonumber \\
&\leq & \pam\left[\lo\left(\tmo,n(1-\zeta_n)\right) < \eta_1'
\lo(\F_{p}(n),(1-\zeta_{n})n) \right] \label{8eq112} .
\end{eqnarray}
To estimate the right hand side of \ref{8eq112}, first we prove
the following lemma
\begin{Lem} \label{3lem1.5.8} For all $0< \xi \leq 1$ there
exists $n_4 \geq n_1$ such that for all $n>n_4$ and $\alpha \in
G_n'$
\begin{eqnarray*}
\pa_{\tmo}\left[\lo\left(\tmo,\ \xi n\right) \geq \frac{n}{(\log
n)^7} \right] \geq 1- \frac{16(\log_p n)^3  }{ \eto^2 \sigma^6
(\log n)^{\frac{1}{2}}\xi }  .
\end{eqnarray*}
\end{Lem}
\begin{Pre}
Let us define the two points $\tM_{<} \in [{M_n}',\tmo] $ and
$\tM_{>} \in [\tmo,{M_n}]$ by
\begin{eqnarray}
& & \tM_<=\sup\left\{t,\ 0>t>{M_n}',\ S_{t}-S_{\tmo} \geq \log n-\left(6+1/2\right)\log_2 n \right\} , \label{10eq79} \\
& &\tM_>=\inf\left\{t,\ 0<t<{M_n},\ S_{t}-S_{\tmo} \geq \log
n-\left(6+1/2\right)\log_2 n \right\}. \label{10eq80}
\end{eqnarray}
Using \ref{hyp4}, \ref{3eq319} and \ref{3eq320} it is easy to show
that for all $n>n_1$ and all $\alpha \in G_n'$ these two points
exist. \\
By the Markov inequality and using \ref{c4s1l1}, we obtain that
\begin{eqnarray*}
& & \pam\left[T_{\tM_<-1} \wedge T_{\tM_>+1} >\xi n \right] \leq
|\tM_<-\tM_>|^3  \exp\left[(S_{\tM_{<}}- S_{\tmo})\vee
(S_{\tM_{>}}- S_{\tmo}) \right] ( \xi n)^{-1} .
\end{eqnarray*}
Using \ref{10eq79} and \ref{10eq80} and the property
\ref{interminibb} we get that for all $n>n_1$ and all $\alpha \in
G_n'$
\begin{eqnarray*}
& & \pam\left[T_{\tM_<-1} \wedge T_{\tM_>+1}>\xi n \right] \leq
\frac{8(\log_p n)^3 }{ \eto^2 \sigma^6 (\log n)^{\frac{1}{2}}\xi
}.
\end{eqnarray*}
Therefore, for all  $n>n_1$ and all $\alpha \in
G_n'$
\begin{eqnarray}
& & \pam\left[\lo\left(\tmo,\xi n \right) \geq \frac{n}{ (\log
n)^7 } \right] \geq
 \pam\left[\lo\left(\tmo,T_{\tM_<-1} \wedge
T_{\tM_>+1}\right) \geq \frac{n}{(\log n)^{7}}
\right]-\frac{8(\log_p n)^3 }{ \eto^2 \sigma^6 (\log
n)^{\frac{1}{2}}\xi } . \label{3eq1.291}
\end{eqnarray}
By Lemma \ref{4.9}
\begin{eqnarray}
& & \pam\left[\lo\left(\tmo,T_{\tM_<-1} \wedge
T_{\tM_>+1}\right) \geq \frac{n}{(\log n)^{7}} \right] \nonumber  \\
& = & \left(1-\alpha_{\tmo}\pa_{\tmo+1}\left[T_{\tmo} \geq
T_{\tM_>+1}\right] -\beta_{\tmo}\pa_{\tmo-1}\left[T_{\tmo} \geq
T_{\tM_<-1}\right]\right)^{\left[n(\log n)^{-7}\right]+1} .
\label{10eq87}
\end{eqnarray}
Using \ref{k1}, \ref{k2} and the definition of $\tM_>$ and $\tM_<$
we have
\begin{eqnarray}
1-\alpha_{\tmo}\pa_{\tmo+1}\left[T_{\tmo}>T_{\tM_>}\right]
-\beta_{\tmo}\pa_{\tmo-1}\left[T_{\tmo}>T_{\tM_<}\right] & \geq
&1- \frac{(\log n)^{6+\frac{1}{2}}}{n} \label{10eq95} .
\end{eqnarray}
Now replacing \ref{10eq95} in \ref{10eq87} and noticing that
$(1-x)^m \geq 1-mx$ for all  $0 \leq x \leq 1$ and $m\geq 1 $ we
have
\begin{eqnarray}
 \pam\left[\lo\left(\tmo,T_{\tM_<} \wedge
T_{\tM_>}\right) \geq \frac{n}{(\log n)^{7}} \right] & \geq &
1-\frac{1}{(\log n)^{1/2}}, \label{10eq96}
\end{eqnarray}
inserting \ref{10eq96} in \ref{3eq1.291} we get the lemma.
\end{Pre}
\noindent \\ Coming back to \ref{8eq112}, using Lemma
\ref{3lem1.5.8} with $\xi=1-\zeta_n$, for all $n>n_4$ and all
$\alpha \in G_n'$
\begin{eqnarray}
& & \pam\left[\lo\left(\tmo,n(1-\zeta_n)\right) < \eta_1'
\lo(\F_{p}(n),(1-\zeta_{n})n)\right] \nonumber \\
& \leq & \pam\left[\lo\left(\tmo,n(1-\zeta_n)\right) < \eta_1'
\lo(\F_{p}(n),(1-\zeta_{n})n),\ \lo\left(\tmo,n(1-\zeta_n)\right) \geq
\frac{n}{ (\log n)^7}\right] \nonumber \\
&+&\frac{16(\log_p n)^3 }{ \eto^2 \sigma^6 (\log
n)^{\frac{1}{2}}(1-\zeta_n) } . \label{8eq113}
\end{eqnarray}
 Let us denote $N_1=\left[ n (\log n)^{-7}\right]+1$, we have
\begin{eqnarray}
& & \pam\left[\lo\left(\tmo,n(1-\zeta_n)\right) < \eta_1'
\lo(\F_{p}(n),(1-\zeta_{n})n),\
\lo\left(\tmo,n(1-\zeta_n)\right) \geq \frac{n}{(\log n)^7}\right] \nonumber \\
& \leq &  \sum_{l=N_1}^n \pam\left[ \lo(\F_{p}(n),\ T_{\tmo,l+1})
> \frac{l}{\eta_1'} \right] \label{eq2829} ,
\end{eqnarray}
recall $T_{\tmo,l+1}=\inf\left\{k>T_{\tmo,l},\ X_k=\tmo
\right\}$ for all $l \geq 1$  and $T_{\tmo,1}\equiv T_{}{\tmo}$ (see \ref{3.7sec}).\\
To estimate the probability in \ref{eq2829} we want to use
exponential Markov inequality. We need the following lemma which
is easy to prove by elementary computations.
\begin{Lem} \label{lemeta}  Let $(\rho_n,n\in \N)$
be a positive decreasing sequence such that $\lim_{n\rightarrow +
 \infty}\rho_n (\log_2 n)^{2}=0$, there exists $n_5
 \geq n_1$ such that for all  $n>n_5$  and all $\alpha \in G_n'$
\begin{eqnarray}
\frac{1}{\eta_1'}-x_1-\rho_n \geq (\nu_n x_1^2)/4 >0. \nonumber
\end{eqnarray}
where $\eta_1'$ is defined just before \ref{eq450}, $x_1$ and
$\nu_n$ just before \ref{eqeq12b}.
\end{Lem}
\noindent Coming back to \ref{eq2829}, we have
\begin{eqnarray}
& &  \sum_{l=N_1}^n \pam\left[ \lo(\F_{p}(n),\ T_{\tmo,l+1})
> \frac{l}{\eta_1'} \right]  \nonumber \\
& \leq  & \sum_{l=N_1}^n   \pam\left[  \lo(\F_{p}(n),\ T_{\tmo,l+1})-(l+1)x_1
> l\left(\frac{1}{\eta_1'}-x_1\left(1+\frac{1}{N_1} \right) \right)
\right] . \label{4eq145}
\end{eqnarray}
Using Lemma \ref{lemeta}, with $\rho_n=\frac{x_1}{N_1} \equiv
\frac{x_1}{\left[ n / (\log n)^{-7} \right]}$ and property
\ref{8eq37}, for all $n>n_5$ and $\alpha \in G_n'$
\begin{eqnarray*}
\frac{1}{\eta_1'}-x_1\left(1+\frac{1}{N_1}\right) \equiv
\frac{1}{\eta_1'}-x_1-\frac{x_1}{N_1}  \geq (\nu_n x_1^2)/4>0  .
\end{eqnarray*}
So for all $n>n_5$ and all $\alpha \in G_n'$ we use the
exponential Markov inequality to estimate the probability in the
right hand side of \ref{4eq145}. Let $\lambda>0$ for all $n>n_5$
and all $\alpha\in G_n'$
\begin{eqnarray}
& &  \pam\left[  \lo(\F_{p}(n),\ T_{\tmo,l+1})-(l+1)x_1
> l\left(\frac{1}{\eta_1'}-x_1\left(1+\frac{1}{N_1} \right) \right)
\right] \nonumber \\
 & \leq & \Eam\left[\exp\left\{\lambda\left(\lo(\F_{p}(n),\
T_{\tmo,\ l+1})-(l+1)x_1\right)\right\}\right] \exp \left[-\lambda
l\left(\frac{1}{\eta_1'}-x_1\left(1+\frac{1}{N_1}\right)\right)\right]
. \label{4eq147}
\end{eqnarray}
By the strong Markov property we have
\begin{eqnarray}
& &  \Eam\left[\exp\left\{\lambda\left(\lo(\F_{p}(n),\ T_{\tmo
,l+1})-(l+1)x_1\right)\right\}\right] \nonumber \\
&=& \left(\Eam\left[\exp\left\{ \lambda\left(\lo(\F_{p}(n),\
T_{\tmo})-x_1\right)\right\}\right]\right)^{l+1} .\label{8eq123}
\end{eqnarray}
To estimate the laplace transform on the right hand side of
\ref{8eq123} we
  use H\"older inequality  and the results of M. Cs\"org\"o L. Horv\'ath and P.
R\'ev\'esz (see \cite{Revesz} pages 279-280) :  choosing
\begin{eqnarray}
\lambda &=& \frac{\left(u^{+}_{n} \wedge u^{-}_{n}
 \right)^2}
{(|\F_{p}^{-}(n)|+|\F_{p}^{+}(n)|)^{3}} , \label{lambdac2}
\end{eqnarray}
where
\begin{eqnarray*}
u^{+}_{n} &=& \min_{ q\in \F_{p}^{+}(n)}\left( \beta_q
\pa_{q-1}\left[T_q>T_{\tmo}\right] \right),\\
u^{-}_{n} &=& \min_{ q\in \F_{p}^{-}(n)}\left( \alpha_q
\pa_{q+1}\left[T_q>T_{\tmo}\right]\right)  , 
\end{eqnarray*}
 $\F_{p}^{-}(n)$ and $\F_{p}^{+}(n)$ have been defined just before \ref{8eq38} and
$a \wedge b =\min(a,b)$, we get
\begin{eqnarray}
 \Eam\left[\exp\left\{\lambda \left( \lo\left(\F_{p}(n),T_{\tmo}\right)-\Eam\left[ \lo\left(\F_{p}(n),T_{\tmo}\right)\right]\right)\right\}\right]  & \leq &
\exp\left[ \frac{2 \lambda }{|\F_{p}(n)|} \right].
\label{propomomexp1eq1}
\end{eqnarray}
\noindent Now using \ref{4eq147}, for  all $n>n_5 $ all $\alpha\in
G_n'$ and all $ l \geq N_1 $
\begin{eqnarray}
& &  \pam\left[  \lo(\F_{p}(n),\ T_{\tmo,l})-(l+1)x_1
> l\left(\frac{1}{\eta_1'}-x_1\left(1+\frac{1}{N_1} \right) \right)
\right] \nonumber \\
 & \leq & \exp \left[-\lambda
l\left\{\frac{1}{\eta_1'}-\left(x_1+\frac{2}{|\F_{p}(n)|}\right)\left(1+\frac{1}{N_1}\right)\right\}\right]
. \label{4eq147b}
\end{eqnarray}
Inserting \ref{4eq147b} in \ref{4eq145} and using \ref{eq2829}, we
deduce that for all $n>n_5$ and $\alpha\in G_n'$
\begin{eqnarray}
& & \pam\left[\lo\left(\tmo,n(1-\zeta_n)\right) < \eta_1'
\lo(\F_{p}(n),(1-\zeta_{n})n),\
\lo\left(\tmo,n(1-\zeta_n)\right) \geq \frac{n}{ (\log n)^7}\right] \nonumber \\
& \leq &
\frac{\exp\left[-N_1\lambda\left(\frac{1}{\eta_1'}-\left(x_1+\frac{2}{|\F_{p}(n)|}\right)\left(1+1/N_1\right)\right)
\right]}{1-\exp\left[-\lambda
\left(\frac{1}{\eta_1'}-\left(x_1+\frac{2}{|\F_{p}(n)|}\right)\left(1+1/N_1\right)\right)\right]}
. \label{10eq137}
\end{eqnarray}
Using Lemma \ref{lemeta} with $\rho_n= x_1/N_1 + (2/|\F_{p}(n)|)
\left(1+1/N_1\right)$ we have for all $n>n_5$ and all $\alpha \in
G_n'$
\begin{eqnarray}
   \frac{1}{\eta_1'}-\left(x_1+\frac{2}{|\F_{p}(n)|}\right)\left(1+1/N_1\right)   \geq  (\nu_n
x_1 ^2)/4 \label{10eq142} .
\end{eqnarray}
Using \ref{10eq142} and \ref{10eq137} we obtain, after an easy
computation, that for all $n>n_5$ and all $\alpha \in G_n'$
\begin{eqnarray}
& & \pam\left[\lo\left(\tmo,n(1-\zeta_n)\right) < \eta_1'
\lo(\F_{p}(n),(1-\zeta_{n})n),\
\lo\left(\tmo,n(1-\zeta_n)\right) \geq \frac{n}{ (\log n)^7}\right] \nonumber \\
& \leq &  \frac{8 \exp\left(-(N_1\lambda   \nu_n x_1^2)/4
\right)}{ \lambda \nu_n x_1^2} . \label{10eq150}
\end{eqnarray}
Now we need a lower and upper bound for $\lambda$, using
\ref{lambdac2} and the properties \ref{8eq38} and \ref{8eq39} we
deduce that for all $n> n_1$ and all $\alpha \in G_n'$ we have
\begin{eqnarray}
  \frac{1}{(g_1(n))^2  (\log_2 n \log_p n)^6} \leq \lambda \leq \frac{2}{  (\log_2 n \log_p n)^6}
, \label{10eq151}
\end{eqnarray}
with  $g_1(n)=\exp\left[( (4 \sqrt{3} \sigma f_{p}(n))^2  \log_3 (n)
)^{1/2} \right]$. We deduce that there exists $n_6\geq n_5 $ such
that for all $n\geq n_6$ and all $\alpha \in G_n'$
\begin{eqnarray}
& & \pam\left[\lo\left(\tmo,n(1-\zeta_n)\right) < \eta_1'
\lo(\F_{p}(n),(1-\zeta_{n})n),\ \lo\left(\tmo,n(1-\zeta_n)\right) \geq
\frac{n}{ (\log n)^7}\right]  \leq \exp(-n^{1/2}) \label{10eq159}
\end{eqnarray}
Collecting \ref{10eq159}, \ref{8eq113}, \ref{8eq112},
\ref{8eq109}, \ref{10eq104} and taking $p=2$, $n_1'=n_2 \vee n_3
\vee n_4 \vee n_6 $ we get \ref{5.31}.
\end{Pre}

\subsection{Proof of Theorem \ref{tpslocmo1nn} }

Clearly the probability in \ref{extension} is bounded from above
by :
\begin{eqnarray*}
 |\mathbb{L}(n)| \max_{ k \in \mathbb{L}(n)} \left\{ \pao\left[\left|\frac{\lo(k,n)}{n}-
\frac{1}{\Ea_{k}\left[\lo(W_n,T_{k})\right]} \right|> \frac{ 1 }{
(\log_2 n)^{1+\rho}}
 \right]\right\}.
\end{eqnarray*}
We are left to give an upper bound of the probability into the
bracket, uniformly in $k \in \F_{p}(n)$. Since the method we use is
exactly the same as the one of the proof of Theorem
\ref{tpslocmo}, we only sketch it. First replace all the $\tmo$ by
$k$ except in "$\F_{p}(n)$", "$\F_{p}^{+}(n)$", "$\F_{p}^{-}(n)$" and "$W_n$" that do
not change. This will change of course the
definitions of $x_1$, $x_2$ and $\lambda$. \\
Then the modifications needed are only based on the following
fact, let $k \in \F_{p}(n)$ and $l \in W_n, l \neq k$, with Q
probability 1 we have
\begin{eqnarray}
 \exp(S_l-S_k) & \equiv & \exp(S_k-S_{\tmo}) \exp(S_l-S_{\tmo}) \nonumber \\
 & \leq  & \left(\log \frac{1-\eta_0}{ \eta_0}\right)^{|k-\tmo|} \exp(S_l-S_{\tmo})   \leq (\log_2 n)
 \exp(S_l-S_{\tmo}).\label{4.80}
\end{eqnarray}
As we will see in section \ref{5.5b}, it is from \ref{4.80} and \ref{8eq36} (respectively \ref{8eq37}) that we get \ref{lok} (respectively \ref{lokk}). \\
The main steps of the proof of Theorem \ref{tpslocmo1} have to be modified as follow : \\
First remark that, using \ref{lok}, \ref{eqeq1} remains true. In
Proposition \ref{8eq33}, \ref{5.31} become  :
\begin{eqnarray*} & & \pa_k\left[\lo(k,n)<\eta_1 n \right]\leq
\frac{4}{\eto} \frac{ 1}{ (\log_p n)^{\rho}}.
\end{eqnarray*}
Notice that, comparing with \ref{5.31} the $\log_2 n$ disappears,
this comes from the fact that we use \ref{lok} (instead of
\ref{2.7}) to prove the equivalent of \ref{4.40} (replacing
$\pa_{\tmo}$ by $\pa_k$). Moreover it is important to notice that
Lemma \ref{lemeta}
remains true and $\lambda$ still verify \ref{10eq151}. \\
Moreover \ref{5.32} became
\begin{eqnarray*}
& & \pa_k\left[\lo(k,n)>\eta_2 n \right]\leq \frac{32(\log_p n)^3
\log_2 n }{ \eto^2 \sigma^6 (\log n)^{\frac{1}{2}} }.
\end{eqnarray*}
since we use \ref{4.80} for the proof of the equivalent of Lemma
\ref{3lem1.5.8} (replacing, as always, $\tmo$ by $k$).
$\blacksquare$

\section{Proof of the good properties for the environment \label{sec4}}

Here we give the main ideas for the proof of Propositions \ref{propoquidonne}, \ref{propoquidonnebbb},  \ref{pro3.10}, \ref{3pro61b} and Proposition 
\ref{profondab}. This section is organized as follow : in section
\ref{ER1} we recall elementary results on sums of independent
random variables, in section \ref{2par22} we prove Lemma \ref{moexiste} and give
standard results on the basic valley, in section \ref{5.3} (respectively \ref{5.2}) we give some details of the proof of Proposition
\ref{3pro61b} (respectively  Proposition \ref{propoquidonneb}).
 Collecting the results of Section \ref{2par22},
Lemma \ref{moexiste}, Propositions \ref{propoquidonneb} and
\ref{3pro61b} we get Proposition \ref{profondab}.  In Section \ref{5.5b} we
sketch the proof of \ref{lok} and \ref{lokk}. In all this
section we assume that the parameter $\gamma$ (see Definition
\ref{thdefval1b} ) is strictly positive.

\subsection{Elementary results \label{ER1}}

We will always work on the right hand side of the origin, that
means with $(S_m,m\in \N)$, by symmetry we obtain the same result
for $m \in \Z_-$.

 \noindent \\ We introduce the following stopping times, for
$a>0$,
\begin{eqnarray}
& & V^+_a\equiv V^+_a(S_j,j\in \N) = \left\{\begin{array}{l}  \inf  \{m \in \N^*,\ S_m \geq a \} , \label{def1}  \\
 + \infty \textrm{, if such a } m \textrm{ does not exist.} \end{array} \right. \label{4.3} \\
& & V^-_a\equiv V^-_a(S_j,j\in \N)= \left\{\begin{array}{l}  \inf  \{m \in \N^*,\ S_m \leq -a \} ,   \\
 + \infty \textrm{, if such a } m \textrm{ does not exist.}
\end{array} \right.
\end{eqnarray}
The following lemma is an immediate consequence of the Wald
equality (see \cite{Neveu})
\begin{Lem} \label{lem101b}  Assume \ref{hyp1bb}, \ref{hyp0} and \ref{hyp4}, let $a>0$, $d>0$ we have
\begin{eqnarray}
& & Q\left[V_{a}^{-} < V_{d}^{+} \right] \leq
\frac{d+\Ie}{d+a+\Ie} , \label{lem101eq2b} \\
& & Q\left[V_{a}^{-} > V_{d}^{+} \right] \leq
\frac{a+\Ie}{d+a+\Ie} , \label{lem101eq3b}
\end{eqnarray}
recall $ \Ie=\log ((1-\eto)(\eto)^{-1})$.
\end{Lem}

\noindent The following lemma is easy to prove when the
$\epsilon_i=\pm 1$ with a probability $1/2$ and is a simple
extension in our more general case

\begin{Lem} \label{maldita} Assume \ref{hyp1bb}, \ref{hyp0} and \ref{hyp4}
hold, there exists $c_{0}\equiv C_{0}(Q)>0$ and $n_0 \equiv
n_0\left(Q\right)$ such that for all $n>n_0$
\begin{eqnarray}
Q\left[V^-_0>r(n) \right] \leq \frac{c_{0}}{\sqrt{r(n)}}
\label{8lem29},
\end{eqnarray}
$(r(n),n)$ is a strictly positive increasing sequence.
\end{Lem}

\subsection{Standard results on the basic valley $\{{M_n},\tmo,{M_n}'\}$
\label{2par22}}

First let us give the main ideas of the proof of  Lemma \ref{moexiste}

\begin{Prele}{\ref{moexiste}}
To prove that $\{M_{n}',m_{n},M_{n}\}\neq \varnothing$ in probability, it is enough to find a valley $\{M',m,M\}$ that satisfies the three properties of Definition \ref{thdefval1b}.  It is easy to show that $\{M'=\bar{V}^{+}_{\Gamma_{n}},m=\bar{m},M=V^{+}_{\Gamma_{n}}\}$ with $\bar{V}^{+}_{\Gamma_{n}}=\sup\{k<0,\ S_{k}>\Gamma_{n}\}$ and $\bar{m}=\inf\left\{|k|>0, S_{k}=\min_{\{\bar{V}^{+}_{\Gamma_{n}} \leq m \leq V^{+}_{\Gamma_{n}}\}}  S_{m}\right \}$ satisfies these properties in probability. Indeed by definition $\{\bar{V}^{+}_{\Gamma_{n}},\bar{m},V^{+}_{\Gamma_{n}}\}$ satisfies  the two first properties of Definition \ref{thdefval1b}. For the third one, assume for simplicity that $\bar{m}>0$, by definition of $\bar{V}^{+}_{\Gamma_{n}}$ and hypothesis \ref{hyp4} we have $ \Gamma_{n} \leq S_{\bar{V}^{+}_{\Gamma_{n}}} \leq \Gamma_{n}+\Ie$ with probability 1. So  we are left to prove that there exists $c_{0}>0$ and $n_{0}\equiv n_{0}(Q)$ such that for all $n\geq n_{0}$
\begin{eqnarray}
Q\left[S_{\bar{V}^{+}_{\Gamma_{n}}}-\max_{0 \leq t \leq
 \bar{m} }  \left(S_t \right) \leq \gamma \log_2 n 
 \right] & \equiv & Q\left[ \Gamma_{n}-\gamma \log_{2} n \leq \max_{0 \leq t \leq
 \bar{m} }  \left(S_t \right) \leq \Gamma_{n}+\Ie \right] \label{redmb}  \nonumber \\
& \leq &
\frac{ c_{0} \gamma \log_2 n}{\log n}. \label{redm}
\end{eqnarray}
To get this upper bound we make the following remark, the event $\left\{\Gamma_{n}-\gamma \log_{2} n \leq \max_{0 \leq t \leq
 \bar{m} }  \left(S_t \right) \right.$ $\left. \leq \Gamma_{n}+\Ie \right\}$ asks to the walk to reach a point larger than $ \Gamma_{n}-\gamma \log_{2} n$ and then to touch a point $S_{\bar{m}}\leq 0$ before a point larger or equal to $\Gamma_{n}+\Ie$. So the probability on the right hand side of \ref{redmb} can be bounded from above by a constant times the probability $Q_{\Gamma_{n}-\gamma \log_{2}n}\left[V^{-}_{0}<V_{\Gamma_{n}+\Ie}\right]$ (where $Q_{y}[\cdots]\equiv Q[\cdots|S_{0}=y] $) which gives \ref{redm} by Lemma \ref{lem101b}. For more details of this computation see for example \cite{Pierreth} pages 56-58.
\end{Prele}

\begin{Lem} \label{8eq18} There exists $c>0$ such that if \ref{hyp1bb}, \ref{hyp0} and
\ref{hyp4} hold, for all $p\geq 2$ there exists $n_0 \equiv
n_0\left(\Ie,\sigma, \E_Q\left[|\epsilon_0|^3\right] \right)$ such
that for all $n>n_0$
\begin{eqnarray}
& & Q\left[{M_n} \leq ( \sigma^{-1} \log n)^2 \log_p n,\ {M_n}'
\geq -( \sigma^{-1} \log n)^2 \log_p n\right] \geq 1- c
R_{p}(n)
 , \label{2eq126} \\
 & & Q\left[{M_n} >\tmo +f_{p}(n) \right] \geq
1-c R_{p}(n)
 \label{eq2623b} , \\
& &  Q\left[ S_{\tM_{1}}-S_{\tm_1} \leq \log n-\gamma \log_2 n,\
S_{\tM_{1}'}-S_{\tm_1'} \leq \log n-\gamma \log_2 n \right] \geq
1- cR_{p}(n).
\end{eqnarray}
See just before \ref{3eq326} for the definitions of $\tM_1$,
$\tm_1$, $\tM_1'$ and $\tm_1'$, $f_{p}(n)$ is given by \ref{F}.
\end{Lem}
\begin{Pre}
\noindent The proof of this lemma is easy and is omitted.
\end{Pre}

\begin{Lem} \label{10eq42} There exists $c>0$ such that if \ref{hyp1bb}, \ref{hyp0} and \ref{hyp4} hold, for all  $p\geq 2$ there exists $n_0 \equiv
n_0(Q)$ such
that for all $n>n_0$
\begin{eqnarray}
& & Q\left[\min_{k\in \F_{p}^{+}(n)}\left(\beta_k
\pa_{k-1}\left[T^{k-1}_{k}>T^{k-1}_{\tmo} \right]\right) \leq
\frac{1}{g_1(n)},\ \tmo>0\right] \leq  c R_{p}(n)  ,\label{8eq42} \nonumber \\
& & Q\left[ \min_{k\in \F_{p}^{-}(n)} \left(\alpha_k
\pa_{k+1}\left[T^{k+1}_{k}>T^{k+1}_{\tmo} \right] \right) \leq
\frac{1}{g_1(n)},\ \tmo>0\right] \leq  c R_{p}(n), \end{eqnarray}
with  $g_1(n)=\exp \left[ \left(
 (2 \sqrt{3} \sigma f_{p}(n) )^2 \log_3 (n) \right)^{1/2} \right]$, recall that  $\F_{p}^{-}(n)$ and $\F_{p}^{+}(n)$ have been defined just before \ref{8eq38} .
\end{Lem}

\begin{Pre}
The proof is omitted it makes use of Lemma \ref{8eq26b} and
elementary facts on sums of i.i.d. random variables.
\end{Pre}

\subsection{Proof of Proposition \ref{propoquidonnebbb}, \ref{Pro}, \ref{pro3.10} and  \ref{3pro61b} \label{5.3} }

\subsubsection{Preliminaries  \label{9.1.2}} \noindent

\noindent
By linearity of the expectation we have :
\begin{eqnarray}
\qquad \Ea_{\tmo}\left[ \lo(W_n,T_{\tmo})\right] &
 \equiv & \sum_{j=\tmo+1}^{M_n}\Ea_{\tmo}\left[ \lo(j,T_{\tmo})\right] +
\sum_{j=M_n'}^{\tmo-1}\Ea_{\tmo}\left[ \lo(j,T_{\tmo})\right] +1 \label{eq23},
\end{eqnarray} 
so using \ref{3.7bb} we get Proposition \ref{Pro}. Now using
Lemma \ref{3.7bb} and hypothesis \ref{hyp4} we easily get the
following lemma
\begin{Lem} \label{9eq13}  Assume \ref{hyp4}, for all  ${M_n}' \leq k \leq {M_n}$, $k\neq \tmo$
\begin{eqnarray}
\frac{\eto}{1-\eto}\frac{1}{e^{S_k-S_{\tmo}}} \leq
\Ea_{\tmo}\left[ \lo(k,T_{\tmo})\right] \leq
\frac{1}{\eto}\frac{1} {e^{S_k-S_{\tmo}} } ,
\end{eqnarray}
with a $Q$ probability equal to one.
\end{Lem}
\noindent Proposition \ref{pro3.10} is a trivial consequence of Lemma \ref{9eq13} and \ref{eq23}. The following lemma is easy
to prove :
\begin{Lem} \label{3lem62} For all $\alpha \in \Omega_1$ and $n>3$, with a $Q$ probability equal to one we have
\begin{eqnarray}
& & \sum_{j=\tmo+1}^{M_n}\frac{1}{e^{S_{j}-S_{\tmo}}} \leq
\sum_{i=1}^{N_{n}+1} \frac{1}{e^{a(i-1)}}\sum_{j=m_n+1}^{M_n}
\un_{S_{j}-S_{\tmo} \in [a(i-1),ai[} ,\label{3eqlem1.209} \\
& & \sum_{j=M_n'}^{m_n-1} \frac{1}{e^{S_{j}-S_{\tmo}}} \leq
\sum_{i=1}^{N_{n}+1}
\frac{1}{e^{a(i-1)}}\sum_{j=M_n'}^{m_n-1} \un_{S_{j}-S_{\tmo} \in
[a(i-1),ai[} ,\label{3eqlem1.210}
\end{eqnarray}
where $a=\frac{\Ie}{4}$, $N_{n}=[(\Gamma_n+\Ie)/a]$, recall that $\Ie=\log ((1-\eto)(\eto)^{-1})$ and $\un$ is the indicator function.
\end{Lem}

\noindent \\ Using Proposition \ref{pro3.10} and Lemma \ref{3lem62}, we have for
all $n>3$
\begin{eqnarray}
 \E_Q\left[ \Ea_{\tmo}\left[ \lo(W_n,T_{\tmo})\right] \right]
\leq  1 &+& \sum_{i=1}^{N_{n}+1} \frac{1}{e^{a(i-1)}}
\E_Q\left[\sum_{j=m_n+1}^{M_n} \un_{S_{j}-S_{\tmo} \in
[a(i-1),ai[} \right] \label{9.161} \\ &+&
\sum_{i=1}^{N_{n}+1} \frac{1}{e^{a(i-1)}}\E_Q\left[
\sum_{j=M_n'}^{m_n-1} \un_{S_{j}-S_{\tmo} \in [a(i-1),ai[}\right]
. \nonumber
\end{eqnarray}

\noindent The next step for the proof of Proposition \ref{3pro61b} is to show that the two expectations $\E_Q[...]$ on
the right hand side of \ref{9.161} are bounded by a constant
depending only on the distribution $Q$ times a polynomial in $i$. This result is given by Lemma \ref{5.7}.

\subsubsection{Proof of Proposition \ref{3pro61b}}

\begin{Rem}
We give some details of the proof of Proposition \ref{proMomo}
mainly because it is based on a very nice cancellation that occurs
between two $\Gamma_{n} \equiv \log n+\gamma \log_{2}n$, see formulas \ref{5.29} and \ref{5.69}.
Similar cancellation is already present in \cite{Kesten2}.
\end{Rem}

\noindent Let us define the following stopping times, let $i>1$ :
\begin{eqnarray*}
& & u_0=0, \\
& & u_1 \equiv V_0^-= \inf\{m>0,\ S_m<0\}, \\
& & u_i = \inf\{m>u_{i-1},\ S_m<S_{u_{i-1}}\}.
\end{eqnarray*}

\noindent \\ The following lemma give a way to characterize the
 point $\tmo$, it is inspired by the work of \cite{Kesten2} and is
 just inspection

\begin{Lem} \label{8eq26b} Let $n>3$ and $\gamma>0$, recall $\Gamma_n =\log n +\gamma \log_2 n$, assume $\tmo > 0 $, for all  $l\in
\N^*$ we have
\begin{eqnarray}
 \tmo=u_l &  \Rightarrow & \left\{
\begin{array}{l}
\bigcap_{i=0}^{l-1}\left\{ \max_{ u_i \leq j \leq u_{i+1} } (S_i)-S_{u_{i}} < \Gamma_n\right\} \textrm{ and } \\
\max_{ u_l \leq j \leq u_{l+1} } (S_i)-S_{u_{i}} \geq \Gamma_n
 \textrm{ and } \\
{M_n}=V^+_{\Gamma_n,l}
\end{array} \right. \label{2eq121b}
\end{eqnarray}
where
\begin{eqnarray}
V^+_{z,l}\equiv V^+_{z,l}\left(S_j,j\geq 1\right) =
\inf\left(m>u_l,\ S_m-S_{u_l} \geq z  \right).
\end{eqnarray}
\end{Lem}

\noindent A similar characterization of $\tmo$ if $\tmo\leq 0$ can
be done (the case $\tmo=0$ is trivial).

\begin{Lem} \label{5.7} There exits $c_{0} \equiv c_{0}(Q)$ such that for all $i\geq 1$
and all $n$ :
\begin{eqnarray}
& & \E_Q\left[\sum_{j=m_n+1}^{M_n} \un_{S_{j}-S_{\tmo} \in
[a(i-1),ai[} \right] \leq c_{0} \times i^2 \label{9.161bb}, \\
& & \E_Q\left[ \sum_{j=M_n'}^{m_n-1} \un_{S_{j}-S_{\tmo} \in
[a(i-1),ai[}\right] \leq c_{0} \times i^2 .  \label{9.161bbb}
\end{eqnarray}
\end{Lem}

\begin{Pre}
We will only prove \ref{9.161bb}, we get \ref{9.161bbb} symmetrically 
moreover we assume that $m_n>0$, computations are similar for the
case $m_n \leq 0$. Thinking on the basic definition of the expectation, we 
 need an upper bound for the probability :
\begin{eqnarray*}
Q\left[\sum_{j=m_n+1}^{M_n} \un_{S_{j}-S_{\tmo} \in
[a(i-1),ai[}=k. \right]
\end{eqnarray*}
First we make a partition over the values of $m_n$ and then we
use Lemma \ref{8eq26b}, we get  :
\begin{eqnarray} 
\quad Q\left[\sum_{j=m_n+1}^{M_n} \un_{S_{j}-S_{\tmo} \in
[a(i-1),ai[}=k. \right] & \equiv &  \sum_{l \geq 0}
Q\left[\sum_{j=m_n+1}^{M_n} \un_{S_{j}-S_{\tmo} \in
[a(i-1),ai[}=k,\ \tmo=u_l \right] \nonumber \\
 & \leq & \sum_{l \geq 0} Q\left[ \A_l^+,\max_{ u_l \leq j \leq u_{l+1} } (S_j)-S_{u_{l}}
\geq \Gamma_n,\A_l^-  \right] \label{5.64}
\end{eqnarray}
where
\begin{eqnarray*}
& & \A_l^+=\sum_{s=u_l+1}^{V^+_{\Gamma_n,l}}\un_{\{S_j-S_{u_l} \in
[a(i-1),ai[\}}=k,\ \\
& & \A_l^- =\bigcap_{r=0}^{l-1}\left\{ \max_{ u_r \leq j \leq
u_{r+1} } (S_r)-S_{u_{r}} < \Gamma_n\right\},\ \A_0^-=\Omega_1.
\end{eqnarray*}
for all $l\geq 0$. By the strong Markov property we have :
\begin{eqnarray}
 Q\left[ \A_l^+,\ \max_{ u_l \leq j \leq u_{l+1} } (S_j)-S_{u_{l}}
\geq \Gamma_n,\ \A_l^-  \right]  \leq  Q\left[ \A_0^+,\
V_0^->V^+_{\Gamma_n} \right] Q\left[\A_l^- \right]. \label{5.63}
\end{eqnarray}
The strong Markov property gives also that the sequence $(\max_{
u_r \leq j \leq u_{r+1} } (S_r)-S_{u_{r}} < \Gamma_n ,r \geq 1)$
 is independent and identically distributed, therefore :
\begin{eqnarray}
Q\left[\A_l^-
\right] \leq \left(Q\left[V_0^-<V^+_{\Gamma_n}\right]\right)^{l-1}. \label{5.60}
\end{eqnarray}
We notice that $Q\left[ \A_0^+,\ V_0^->V^+_{\Gamma_n} \right]$
 does not depend on $l$, therefore, using \ref{5.64}, \ref{5.63}
 and \ref{5.60} we get :
\begin{eqnarray}
\quad Q\left[\sum_{j=m_n+1}^{M_n} \un_{S_{j}-S_{\tmo} \in
[a(i-1),ai[}=k. \right] & \leq & (1+(Q\left[V_0^-\geq
V^+_{\Gamma_n}\right])^{-1})Q\left[ \A_0^+,\ V_0^->V^+_{\Gamma_n}
\right] \label{5.29}
\end{eqnarray}
To get an upper bound for $Q\left[ \A_0^+,\ V_0^->V^+_{\Gamma_n}
\right]$, first we introduce the following sequence of stopping
times, let $s>0$ :
\begin{eqnarray*}
& & H_{ia,0}=0, \\
& & H_{ia,s}=\inf\{m>H_{ia,s},\ S_m \in [(i-1)a,ia[\}.
\end{eqnarray*}
Making a partition over the values of $H_{ia,k}$, by the Markov property we get:
\begin{eqnarray}
& & Q\left[ \A_0^+,\ V_0^->V^+_{\Gamma_n} \right] \nonumber \\
& \leq & \sum_{w \geq 0}\int_{(i-1)a}^{ia} Q\left[H_{ia,k}=w,S_w \in
dx,\bigcap_{s=0}^{w}\{S_s>0\},\bigcap_{s=w+1}^{\inf \{l>w,S_l
\geq \Gamma_n\}}\{S_s>0\} \right]  \nonumber \\
& \leq &  Q\left[H_{ia,k}<V^-_0\right] \max_{ (i-1)a \leq x \leq
ia } \left\{Q_x\left[V_{\Gamma_n-x}^{+}<V^-_{x} \right]\right\} \equiv
Q\left[H_{ia,k}<V^-_0\right] Q_{ia}\left[V_{\Gamma_n-ia}^{+}<V^-_{ia}
\right] \label{5.69}
\end{eqnarray}
To finish we need an upper bound for
$Q\left[H_{ia,k}<V^-_0\right]$, we do not want to give details of
the computations for this because it is not difficult, however the
reader can find these details in \cite{Pierreth} pages 142-145. We
have for all $i>1$:
\begin{eqnarray}
\qquad  & &  Q\left[H_{ia,k}<V^-_0\right]  \leq
Q\left[V^-_0>V^+_{(i-1)a}\right]
\left(1-Q\left[\epsilon_0<-\frac{\Ie}{2}\right]Q_{(i-1)a-\frac{\Ie}{4}}\left[V^+_{(i-1)a}
\geq V^-_0\right] \right)^{k-1} \ ,  \label{5.70} \\
 &  &
Q\left[H_{a,k}<V^-_0\right]  \leq Q\left[\epsilon_0\geq 0\right]
\left(1-Q\left[\epsilon_0<-\frac{\Ie}{4}\right] \right)^{k-1} \ .  \label{5.68}
\end{eqnarray}
So using \ref{5.29}, \ref{5.69}-\ref{5.68}, and Lemma \ref{lem101b} one can find a constant $c_0$ that
depends only on the distribution $Q$ such that for all $i\geq 0$ : 
\begin{eqnarray*}
& & \E_Q\left[\sum_{j=m_n+1}^{M_n} \un_{S_{j}-S_{\tmo} \in
[a(i-1),ai[} \right]  \equiv  \sum_{k=1}^{+ \inf } k
Q\left[\sum_{j=m_n+1}^{M_n} \un_{S_{j}-S_{\tmo} \in
[a(i-1),ai[}=k \right]    \leq c_0 \times  i^2,
\end{eqnarray*}
which provide \ref{9.161bb}
\end{Pre}

\subsubsection{Proof of Proposition \ref{propoquidonnebbb}}

To show Proposition \ref{propoquidonnebbb} we use the same method previously used to prove Proposition \ref{3pro61b}, the key point is to show 
 the following Lemma :  

\begin{Lem} There exits a constant $c_{0}\equiv c_{0}(Q)$ such that for all strictly positive increasing sequences $(\theta(n),n)$, $i\geq 1$
and all $n$ :
\begin{eqnarray}
& & \E_Q\left[\sum_{j=m_n+\theta(n)}^{M_n} \un_{S_{j}-S_{\tmo} \in
[a(i-1),ai[} \right] \leq \frac{c_{0} \times i^2}{\sqrt{\theta(n)}} \label{9.161bb1}, \\
& & \E_Q\left[ \sum_{j=M_n'}^{m_n-\theta(n)} \un_{S_{j}-S_{\tmo} \in
[a(i-1),ai[}\right] \leq \frac{c_{0} \times i^2}{\sqrt{\theta(n)}}.  \label{9.161bbb2}
\end{eqnarray}
\end{Lem}
\noindent
We will not give the details of this proof because it is very similar to the proof of Lemma \ref{propoquidonneb}.
Just notice that $\sqrt{\theta(n)}$ comes from the fact that the probability $Q
\left[\bigcap_{s=1}^{\theta_{n}}\{S_{s}>0\}\right]\equiv Q[V_{0}>\theta(n)] $ appears when we give an upper bound of $Q\left[\sum_{j=m_n+\theta(n)}^{M_n} \un_{S_{j}-S_{\tmo} \in
[a(i-1),ai[}=k\right]$. This comes from the fact that the event $\bigcap_{l=1}^{\theta(n)}\{S_{\tmo+l}-S_{\tmo}>0\}$ is hidden in the Definition of $m_{n}$.
A last remark, the upper bound we get here is good for sequences $(\theta(n),n)$ that grow very slowly. In the next section we use another method more powerful for sequences that grow more rapidly.

\subsection{Proof of Proposition \ref{propoquidonneb} \label{5.2}}

\subsubsection{Preliminaries  \label{9.1.2}} \noindent
\noindent
It is here that the explicit form of $f_{p}(n)$ given by \ref{F} will
become clear. Using \ref{eq23} and Lemma \ref{9eq13} we need only 
to find an upper bound for
\begin{eqnarray}
\sum_{l=\tmo+f_{p}(n)+1}^{{M_n}}\frac{1}{e^{S_l-S_{\tmo}}}+\sum_{l=M_n'}^{\tmo-f_{p}(n)-1}\frac{1}{e^{S_l-S_{\tmo}}}
. \label{lessomcv}
\end{eqnarray}
Assume for the moment that
\begin{eqnarray}
S_k-S_{\tmo} \geq 2\log(|k-\tmo|),\ \forall k \in
\{{M_n}',\cdots,\tmo-f_{p}(n),\ \tmo+f_{p}(n),\cdots,{M_n}'\}  ,
\label{cvtermpr}
\end{eqnarray}
with a $Q$ probability larger than $1-cR_{p}(n)$ with $c>0$. Then for all  $k \in
\{{M_n}',\cdots,\tmo-f_{p}(n),\ \tmo+f_{p}(n),\cdots,{M_n}'\}$, we have
\begin{eqnarray*}
\frac{1}{e^{S_k-S_{\tmo}}} \leq \frac{1}{(k-\tmo)^2} .
\end{eqnarray*}
This implies the convergence of the two partial sums in
\ref{lessomcv} and therefore with a $Q$ probability close to 1
\begin{eqnarray*}
 \Ea_{\tmo}\left[ \lo(V_n^{c,r},T_{\tmo})\right]  & \leq &  \frac{1}{\eto} \left(
 \sum_{l=\tmo+f_{p}(n)+1}^{{M_n}}\frac{1}{e^{S_l-S_{\tmo}}}+\sum_{l=M_n'}^{\tmo-f_{p}(n)-1}\frac{1}{e^{S_l-S_{\tmo}}}
 \right)
\\ & \leq & \frac{1}{\eto}\frac{2}{f_{p}(n)+1} ,
\end{eqnarray*}
this entails \ref{2.7}. \ref{cvtermpr} is a consequence of
 Propositions \ref{proMomo} proved in the following section.

\subsubsection{Study of the potential $(S_m,\ m\in \Z)$ in $\bar{\F}_{p}(n)$ \label{9.1.3}} \noindent

\noindent We will assume $\tmo > 0$,  the case $\tmo \leq 0$ is
similar.

\begin{Pro} \label{proMomo} There exists $c>0$such that if \ref{hyp1bb}, \ref{hyp0} and \ref{hyp4} hold,  for all $p\geq 2$
 there exists $c_{0} \equiv c_{0}(Q)$ and $n_0 \equiv n_{0}(Q)$
such that for all  $n>n_0$
\begin{eqnarray}
&  & Q\left[\bigcup_{k=\tmo+f_{p}(n)}^{M_n}\left\{ S_{k}-S_{\tmo}
\leq
2\log (k-\tmo) \right\}, \ \tmo>0\right]   \leq c R_{p}(n)+ \frac{c_{0}}{\log_{p}n},  \label{1eq21} \\
&  & Q\left[\bigcup_{k=M_{n'}}^{\tmo+f_{p}(n)}\left\{ S_{k}-S_{\tmo}
\leq
2\log (\tmo-k) \right\}, \ \tmo>0\right]   \leq cR_{p}(n)+ \frac{c_{0}}{\log_{p}n}. \label{1eq22}
\end{eqnarray}
\end{Pro}


\begin{Prepr}{\ref{proMomo}} We will only prove \ref{1eq21} we get \ref{1eq22} with the same arguments. Let $n\geq 3$, we denote $L(n)=(
\sigma^{-1} \log n)^2 \log_p n$ and $[L(n)]$ the integer part of
$[L(n)]$. By \ref{2eq126} and \ref{eq2623b}, there exists $n_1
\equiv n_1\left(\Ie,\sigma,\E_Q\left[|\epsilon_0|^3\right]
\right)$ such that for all $n>n_1$
\begin{eqnarray}
& & Q\left[\bigcup_{k=\tmo+f_{p}(n)}^{M_{n}}\left\{ S_{k}-S_{\tmo}
\leq 2\log (k-\tmo) \right\},\ \tmo>0 \right]  \leq cR_{p}(n)+ 
Q\left[\A \right] , \label{doublesom}
\end{eqnarray}
where
\begin{eqnarray}
 \A= \left\{\bigcup_{k=\tmo+f_{p}(n)}^{M_{n}}\left\{
S_{k}-S_{\tmo} \leq 2\log (k-\tmo) \right\},\ \tmo+f_{p}(n)<M_n \leq
L(n)\right\} . \nonumber
\end{eqnarray}
Using the same method details for the Proof of Lemma \ref{5.7} (from line \ref{5.64} to \ref{5.29}), we get :
\begin{eqnarray}
 Q\left[\A\right] \leq (1+(Q\left[V^{-}_{0} \geq V^{+}_{\Gamma_{n}}\right])^{-1} )Q\left[\A_{0}^{+}\right] \label{eq547}
\end{eqnarray}
where 
\begin{eqnarray}
 \A_{0}^{+}= \left\{\bigcup_{j=f_{p}(n)}^{V^{+}_{\Gamma_{n}}}\left\{
S_{j} \leq 2\log (j) \right\},\ f_{p}(n)< V^{+}_{\Gamma_{n}} \leq
L(n)\right\} . \nonumber
\end{eqnarray}
Let us denote $\Ct_n=\bigcup_{j=f_{p}(n)}^{V^+_{\Gamma_n}}\left\{
S_j \leq 2\log j \right\}$, to estimate $Q\left[\A_{0}^{+}\right]$ we make a partition over the values
of $S_{f_{p}(n)}$,
\begin{eqnarray}
&& Q\left[\A_{0}^{+} \right] \label{eq264511} \nonumber \\
&= & Q\left[\Ct_n,\ V^-_{0}>V^+_{\Gamma_n},\ f_{p}(n)<V^+_{\Gamma_n} \leq L(n) ,
\ S_{f_{p}(n)} \leq 2 \log f_{p}(n) \right] + \label{eq2645}
\\
& & Q\left[\Ct_n,\ V^-_{0}>V^+_{\Gamma_n},\ f_{p}(n)<V^+_{\Gamma_n} \leq L(n),\
S_{f_{p}(n)} > 2 \log f_{p}(n) \right] \label{eq2645bis} .
\end{eqnarray}
 First let us estimate \ref{eq2645}, we remark
that $ \{V^-_{0}>V^+_{\Gamma_n},\  f_{p}(n)<V^+_{\Gamma_n}\leq L(n)\}
\Rightarrow S_{f_{p}(n)} \geq 0 $, so by the Markov property
\begin{eqnarray}
& & Q\left[\Ct_n,\ V^-_{0}>V^+_{\Gamma_n},\ f_{p}(n)<V^+_{\Gamma_n}\leq L(n), \
S_{f_{p}(n)} \leq 2
\log f_{p}(n) \right] \nonumber \\
& \leq & \int_{0}^{2\log
f_{p}(n)}Q\left[ V^-_{0}>V^+_{\Gamma_n},\ f_{p}(n)<V^+_{\Gamma_n}\leq L(n), \ S_{f_{p}(n)} \in dy \right] \nonumber \\
&=& \int_{0}^{2\log f_{p}(n)}Q\left[\bigcap_{k=1}^{f_{p}(n)} \left\{S_k
\geq 0\right\}, \ S_{f_{p}(n)}\in dy \right]
 Q_{y}\left[V^-_0>V^+_{\Gamma_n-y } \right]  \label{4.136}
\end{eqnarray}
where $Q_y[.]=Q[.|S_0=y]$. Since $Q_{y}\left[V^-_0>V^+_{\Gamma_n-y
} \right] $ is increasing en $y$, the term in the right hand side
of \ref{4.136} is bounded by
\begin{eqnarray}
& &  Q_{2 \log f_{p}(n)}\left[V^-_{0}>V^+_{\Gamma_n-2\log (f_{p}(n))
}\right] \int_{0}^{2\log f_{p}(n)}Q\left[\bigcap_{k=1}^{f_{p}(n)}
\left\{S_k \geq 0\right\}, \ S_{f_{p}(n)} \in dy \right] \nonumber
\\
&=& Q\left[V_0^->f_{p}(n)\right]Q_{2 \log
f_{p}(n)}\left[V^-_{0}>V^+_{\Gamma_n-2\log f_{p}(n) }\right]
\label{eq2659} .
\end{eqnarray}
To estimate \ref{eq2645bis}, let us define the stopping time
$U_{f_{p}(n)}=\inf\left\{m > f_{p}(n),\ S_m \leq \log m\right\}$. We
remark that
\begin{eqnarray}
\left\{ \Ct_n, \ S_{f_{p}(n)} \geq 2 \log f_{p}(n)\right\} \Rightarrow
\left\{ f_{p}(n) \leq U_{f_{p}(n)} \leq V^+_{\Gamma_n}\right\} . \nonumber
\end{eqnarray}
Defining
\begin{eqnarray}
 \A'(l) & = & \left\{\Ct_n,\ V^-_{0}>V^+_{\Gamma_n},\ f_{p}(n)
<V^+_{\Gamma_n},\ S_{f_{p}(n)} \geq 2 \log f_{p}(n),\ U_{f_{p}(n)}=l \right\} \nonumber
,
\end{eqnarray}
 we have
\begin{eqnarray}
 Q\left[\Ct_n,\ V^-_{0}>V^-_{\Gamma_n},\ f_{p}(n) <V^+_{\Gamma_n} ,\
S_{f_{p}(n)} \geq 2 \log f_{p}(n) \right] &=&
\sum_{l=f_{p}(n)}^{[L(n)]}Q\left[\A'(l) \right] . \label{9eq53}
\end{eqnarray}
Since by hypothesis  $ Q\left[- \Ie \leq \epsilon_0 \leq \Ie
\right]=1$, we have $U_{f_{p}(n)}=l \Rightarrow  \log l-\Ie \leq S_l
\leq \log l\quad Q.a.s.$, so
\begin{eqnarray}
& & \sum_{l=f_{p}(n)}^{[L(n)]}Q\left[A'(l) \right] \nonumber \\
 &=& \sum_{l=f_{p}(n)}^{[L(n)]}\int_{\log l-\Ie}^{\log
l}Q\left[\Ct_n,\ V^-_{0}>V^+_{\Gamma_n}, f_{p}(n) <V^+_{\Gamma_n},\
S_{f_{p}(n)} \geq 2 \log f_{p}(n),\D_{l,y} \right] . \nonumber
\end{eqnarray}
with $\D_{l,y}=\left\{\bigcap_{j=f_{p}(n)}^{l-1}\left\{S_j < \log j
\right\},\ S_l \in dy \right\} $. By the Markov property, we get
\begin{eqnarray*}
 \sum_{l=f_{p}(n)}^{[L(n)]}Q\left[A'(l) \right]
 &\leq& \sum_{l=f_{p}(n)}^{[L(n)]}\int_{\log l-\Ie }^{\log l} Q_y\left[
V^-_{0}>V^+_{\Gamma_n-y}\right] \times \\
 & &Q\left[ S_{f_{p}(n)} \geq 2 \log
f_{p}(n),\ \bigcap_{j=0}^{l}\left\{ S_j \geq 0 \right\}, \
\D_{l,y}\right].  
\end{eqnarray*}
Using that $Q_y\left[ V^-_{0}>V^+_{\Gamma_n-y}\right]$ is
increasing in $y$, we obtain
\begin{eqnarray}
\sum_{l=f_{p}(n)}^{[L(n)]}Q\left[A'(l) \right] &\leq &
\sum_{l=f_{p}(n)}^{[L(n)]} Q_{\log l}\left[V^-_{0}>V^+_{\Gamma_n-\log
l}\right]   \int_{\log l-\Ie}^{\log l}
Q\left[\bigcap_{i=1}^{l}\left\{S_i>0\right\},
\D_{l,y} \right] \nonumber  \\
& \leq &  Q_{\log
([L(n)])}\left[V^-_{0}>V^+_{\Gamma_n-\log([L(n)])}\right]
\sum_{l=f_{p}(n)}^{[L(n)]}
Q\left[\bigcap_{i=1}^{f_{p}(n)}\left\{S_i>0\right\},\ U_{f_{p}(n)}=l
\right] \nonumber \\
&=&  Q_{\log
([L(n)])}\left[V^-_{0}>V^+_{\Gamma_n-\log([L(n)])}\right]Q\left[V^-_0>f_{p}(n)
\right] \label{eq2665} .
\end{eqnarray}
Using  \ref{eq2665}, \ref{9eq53}, \ref{eq2659} and \ref{4.136} we get  
\begin{eqnarray}
 & & Q\left[\A_{0}^{+}\right]   \nonumber \\
 &\leq & Q\left[V^-_0>f_{p}(n) \right]\left(Q_{2 \log f_{p}(n)}\left[V^-_0>V^+_{\Gamma_n-2\log f_{p}(n)}\right] +  Q_{\log ([L(n)])}
 \left[V^-_0>V^+_{\Gamma_n-\log([L(n)])}\right]\right) \label{3eq64} .
\end{eqnarray}
Now using Lemmata \ref{lem101b} and \ref{maldita}, there exists $n_2 \equiv n_2\left(Q \right)$ and $c_{0} \equiv c_{0}(Q)>0$ such
that  for all $n>n_{2}$ :
\begin{eqnarray}
 & &  Q\left[\A_{0}^{+}\right]
\leq  \frac{ c_{0} \log_{2}(n) }{\Gamma_{n} \sqrt{f_{p}(n)} } \label{canc}
\end{eqnarray}
\noindent Using \ref{canc}, \ref{eq547},  once again Lemma \ref{lem101b}, \ref{doublesom} and finally  
 taking $n_0'= n_1 \vee n_2 $ ends the
 proof of \ref{1eq21}.
\end{Prepr}

\subsection{Proof of \ref{lok} and \ref{lokk} \label{5.5b}}

For \ref{lok}, we use the fact that Lemma \ref{9eq13} remains true
if we replace $\tmo$ by some $l\in \F_{p}(n)$. Then we use \ref{4.80}
and finally we follow the method used in section \ref{5.2}.
\ref{lokk} is obtained in a similar way.

\noindent \\ \textbf{ Acknowledgment : } This article is a part of
my Phd thesis made under the supervision of P. Picco. I would like
to thank him for helpful discussions all along the last three
years. I would like to thank R. Correa and all the members of the
C.M.M. (Santiago, Chili) for their hospitality during all the year
2002. I would like to thank the anonymous referee for his nice comments.

{\small \bibliography{article}}

\vspace{1cm} \noindent
\begin{tabular}{ll}
Centre de Physique Théorique - C.N.R.S.  & \hspace{1cm} Centro de Modelamiento Matem\'atico - C.N.R.S. \\
Université Aix-Marseille II & \hspace{1cm} Universidad de Chile  \\
Luminy Case 907 & \hspace{1cm} Blanco Encalada 2120 piso 7 \\
13288 Marseille cedex 09, France & \hspace{1cm} Santiago de Chile
\end{tabular}

\end{document}